\documentclass[11pt]{article}
\usepackage{mathrsfs}
\usepackage{latexsym, amsfonts,mathrsfs, amsmath, amssymb, amscd, epsfig}

\newcommand{\Div}{\text{div}}
\newtheorem{thm}{Theorem}

\newtheorem{prop}[thm]{Proposition}
\newtheorem{lem}[thm]{Lemma}

\newtheorem{rmk}{Remark}

\newenvironment{pf}{{\noindent \it \bf Proof:}}{{\hfill$\Box$}\\}

\begin{document}

\def\del{\partial}
\def\DOT{\!\cdot\!}
\def\dt{{\Delta t}}
\def\dx{{\Delta x}}
\def\dy{{\Delta y}}
\def\dz{{\Delta z}}
\def\u{{\mbox{\boldmath $u$}}}
\def\g{{\mbox{\boldmath $g$}}}
\def\bu{{\mbox{\boldmath $u$}}}
\def\x{{\mbox{\boldmath $x$}}}
\def\n{{\mbox{\boldmath $n$}}}
\def\bn{{\mbox{\boldmath $n$}}}
\def\veceps{\mbox{\boldmath $\epsilon$}}
\def\vecdel{\mbox{\boldmath $\delta$}}
\def\vecphi{\mbox{\boldmath $\phi$}}
\def\vecpsi{\mbox{\boldmath $\psi$}}
\def\ds{\mbox{d \boldmath $\gamma$}}
\def\vecomi{\mbox{\boldmath $\theta$}}
\def\eps{\varepsilon}
\def\oneb{{1\mbox \tiny b}}
\def\twob{{2\mbox \tiny b}}
\def\disp{\displaystyle}
\def\dspace{\displaystyle \vspace{.15in}}
\def\half{\textstyle \frac12}
\def\Circ{{\footnotesize$\bigcirc$}}
\def\Triangle{{\small$\triangle$}}
\def\Bullet{{\large$\bullet$}}
\def\exac{{\mbox \tiny e}}
\newcommand{\Dtil}{\widetilde{D}}

\def\TIME{\!\times\!}

\title{Existence of Global Steady Subsonic Euler Flows\\
 through Infinitely Long Nozzles}

\author {Chunjing \, Xie\thanks{Department of mathematics, University
of Michigan, 530 Church Street, Ann Arbor, MI 48109-1043 USA.
E-mail: cjxie@umich.edu}\,\, and\,\, Zhouping Xin\thanks{The
Institute of Mathematical Sciences, The Chinese University of Hong
Kong, Shatin, N.T., Hong Kong. E-mail: zpxin@ims.cuhk.edu.hk}}

\date{}
\maketitle

\bigskip

{\bf Abstract:} In this paper, we study the global existence of
steady subsonic Euler flows through infinitely long nozzles without
the assumption of irrotationality. It is shown that when the
variation of Bernoulli's function in the upstream is sufficiently
small and mass flux is in a suitable regime with an upper critical
value, then there exists a unique global subsonic solution in a
suitable class for a general variable nozzle. One of the main
difficulties for the general steady Euler flows, the governing
equations are a mixed elliptic-hyperbolic system even for uniformly
subsonic flows. A key point in our theory is to use a stream
function formulation for compressible Euler equations. By this
formulation, Euler equations are equivalent to a quasilinear second
order equation for a stream function so that the hyperbolicity of
the particle path is already involved. The existence of solution to
the boundary value problem for stream function is obtained with the
help of the estimate for elliptic equation of two variables. The
asymptotic behavior for the stream function is obtained via a blow
up argument and energy estimates. This asymptotic behavior, together
with some refined estimates on the stream function, yields the
consistency of the stream function formulation and thus the original
Euler equations.

\section{Introduction and Main Results}\label{SEIntroduction}

Multidimensional gas flows give rise many outstanding challenging
problems. Since the solutions for the unsteady compressible Euler
equations develop singularities in general\cite{Sideris}, it is not
yet known which function space is suitable to study their
wellposedness\cite{Rauch}. It is natural to start from the steady
Euler equations to understand some important true multidimensional
flow patterns. However, the steady Euler equations themselves are
not easy to tackle, since the equations may not only be hyperbolic
or hyperbolic-elliptic coupled system, but also have discontinuous
solutions such as shock waves and vortex sheets. Therefore, a lot of
approximate models were proposed to study fluid flows. An important
approximate model is the potential flow, which originates from the
study for flows without vorticity. Since 1950's, tremendous progress
has been made on the study for potential flows. Subsonic potential
flows around a body were studied extensively by
Shiffman\cite{Shiffman}, Bers\cite{Bers1,Bers2}, Finn,
Gilbarg\cite{FG1,FG2}, and Dong\cite{Dong}, et al. Subsonic-sonic
flows around a body were established recently by Chen, et al
\cite{CDSW} via compensated compactness method. Significant progress
on transonic flows was made by Morawetz. She first showed the
nonexistence of smooth transonic flows in
general\cite{Morawetz1,Morawetz2,Morawetz3,Morawetz4}, and later
worked on existence of weak solutions to transonic flows by the
theory of compensated compactness\cite{Morawetzcpt1,Morawetzcpt2}.
Existence and stability of transonic shocks in a nozzle for
potential flows were achieved recently with prescribed potential at
downstream in \cite{CF1,CF2}. Xin and Yin obtained existence and
nonexistence of transonic shocks in a bounded nozzle with prescribed
pressure at downstream was obtained in \cite{XY,XY2}. Recently,
well-posedness for subsonic and subsonic-sonic potential flows
through infinitely long 2-D and 3-D axially symmetric nozzles, was
established in \cite{XX1,XX2}. For the study on other aspects on
subsonic potential flows, please refer to \cite{Feistauer,G,GS}.

Besides the potential flow, there is another important approximate
model to compressible Euler equations, incompressible Euler
equations, which approximate to the compressible Euler equations for
flows with small Mach numbers. For the study on the existence of
steady incompressible Euler flows in a bounded domain, please refer
to \cite{Alber,Troshkin,Glass}, etc, and references therein.

For the full compressible Euler equations, well-posedness and
nonexistence of a transonic shock in bounded nozzles with prescribed
pressure at downstream has been obtained in
\cite{XYY,XY3,LXY1,LXY2}. Existence and stability of transonic
shocks in nozzles with prescribed velocity at downstream was shown
in \cite{CCS}.

In this paper, we study the existence of global steady subsonic
Euler flows through  general infinitely long nozzles.

Consider the 2-D steady isentropic Euler equations
\begin{eqnarray}
&&(\rho u)_{x_1}+(\rho v)_{x_2}=0,\label{EEulercontinuityeq}\\
&&(\rho u^2)_{x_1}+(\rho uv)_{x_2}+p_{x_1}=0,\label{EEulermomentumeq1}\\
&&(\rho uv)_{x_1}+(\rho
v^2)_{x_2}+p_{x_2}=0,\label{EEulermomentumeq2}
\end{eqnarray}
where $\rho$, $(u,v)$, and $p=p(\rho)$ denote the density,
velocity and  pressure respectively. In general, it is assumed
that $p'(\rho)>0$ for $\rho>0$ and $p''(\rho)\geq 0$, where
$c(\rho)=\sqrt{p'(\rho)}$ is called the sound speed. The most
important examples include polytropic gases and isothermal gases.
For polytropic gases, $p=A\rho^{\gamma}$ where $A$ is a constant
and $\gamma$ is the adiabatic constant with $\gamma>1$; and for
isothermal gases, $p=c^2\rho$ with constant sound speed $c$
\cite{CF}.

We consider flows through an infinitely long nozzle given by
\begin{equation*}
\Omega=\{(x_1,x_2)|f_1(x_1)<x_2<f_2(x_1),-\infty<x_1<\infty\},
\end{equation*}
which is bounded by
$S_i=\{(x_1,x_2)|x_2=f_i(x_1),-\infty<x_1<\infty\}$, ($i=1,2$).
Suppose that $S_i(i=1,2)$ satisfy
\begin{eqnarray}
&&f_2(x_1)>f_1(x_1)\,\, \text{for}\,\, x_1\in(-\infty,\infty),\label{Enontrivialbdy}\\
&&f_1(x_1)\rightarrow 0,\quad f_2(x_1)\rightarrow 1, \qquad
\text{as}\,\, x_1\rightarrow
-\infty,\label{Eboundary1}\\
&&f_1(x_1)\rightarrow a,\quad f_2(x_1)\rightarrow b>a,\qquad
\text{as}\,\, x_1\rightarrow +\infty,\label{Eboundary2}
\end{eqnarray}
and
\begin{equation}\label{Eboundary3}
\|f_i\|_{C^{2,\alpha}(\mathbb{R})}\leq C\,\,\text{for some}\,\,
\alpha>0\,\,\text{and}\,\, C>0.
\end{equation}
It follows that $\Omega$ satisfies the uniform exterior sphere
condition with some uniform radius $r>0$.

Suppose that the nozzle walls are impermeable solid walls so that
the flow satisfies the no flow boundary condition
\begin{equation}\label{Enoflowbc}
(u,v)\cdot \vec{n}=0\,\,\text{on}\,\,\partial\Omega,
\end{equation}
where $\vec{n}$ is the unit outward normal to the nozzle wall. It
follows from (\ref{EEulercontinuityeq}) and (\ref{Enoflowbc}) that
\begin{equation}\label{Emassflux}
\int_{l}(\rho u,\rho v)\cdot \vec{n} dl \equiv m
\end{equation}
holds for some constant $m$, which is called the mass flux, where
$l$ is any curve transversal to the $x_1-$direction, and $\vec{n}$
is the normal of $l$ in the positive $x_1$-axis direction.

Due to the continuity equation, when the flow is away from the
vacuum, the momentum equations are equivalent to
\begin{eqnarray}
uu_{x_1}+vu_{x_2}+h(\rho)_{x_1}=0\label{Enonconservemeq1},\\
uv_{x_1}+vv_{x_2}+h(\rho)_{x_2}=0\label{Enonconservemeq2},
\end{eqnarray}
where $h(\rho)$ is the enthalpy of the flow satisfying
$h'(\rho)=p'(\rho)/\rho$. So $h(\rho)$ is determined up to a
constant. In this paper, for example, we always choose $h(0)=0$ for
polytropic gases and $h(1)=0$ for isothermal gases. After
determining this integral constant, we denote $B_0=\inf_{\rho>0}
h(\rho)$.

It follows from (\ref{Enonconservemeq1}) and
(\ref{Enonconservemeq2}) that
\begin{equation}\label{EstreamBernoulli}
(u,v)\cdot \nabla (h(\rho)+\frac{1}{2}(u^2+v^2))=0.
\end{equation}
This implies that $\frac{u^2+v^2}{2}+h(\rho)$, which will be called
Bernoulli's function, is a constant along each streamline. For Euler
flows in the nozzle, we assume that in the upstream, Bernoulli's
function is given, i.e.,
\begin{equation}\label{EasymtoticBernoulli}
\frac{u^2+v^2}{2}+h(\rho)\rightarrow B(x_2)\,\, \text{as}\,\,
x_1\rightarrow -\infty,
\end{equation}
where $B(x_2)$ is a function defined on $[0,1]$.

Now let us state our main results in the paper
\begin{thm}\label{EThexistence}
Let the nozzle satisfy (\ref{Enontrivialbdy})-(\ref{Eboundary3}) and
$ \underline{B}>B_0$. There exists a $\delta_0>0$ such that if
\begin{equation}\label{EassumptiononBernoulli}
\inf_{x_2\in[0,1]}B(x_2)=\underline{B},\,\, B'(0)\leq
0,\,\,B'(1)\geq 0\,\,\,\,\text{and}\,\,\|B'(x_2)\|_{C^{0,1}([0,1])}=
\delta\leq \delta_0,
\end{equation}
then there exists $\hat{m}\geq 2\delta_0^{1/8}$ such that for any
$m\in (\delta^{1/4},\hat{m})$,
\begin{enumerate}
\item (Existence) there exists a flow satisfying the Euler
equations (\ref{EEulercontinuityeq})-(\ref{EEulermomentumeq2}),
the boundary condition (\ref{Enoflowbc}),   mass flux condition
(\ref{Emassflux}), and the asymptotic condition
(\ref{EasymtoticBernoulli}); \item (Subsonic flows and positivity
of horizontal velocity)  the flow is globally uniformly subsonic
and has positive horizontal velocity in the whole nozzle, i.e.,
\begin{equation}\label{Esubsonic}
\sup_{\bar{\Omega}}(u^2+v^2-c^2(\rho))<0\,\,\text{and}\,\,
u>0\,\,\text{in}\,\,\bar{\Omega};
\end{equation}
\item (Regularity and far fields behavior)
Furthermore, the flow satisfies
\begin{equation}
\|\rho\|_{C^{1,\alpha}(\Omega)},\|u\|_{C^{1,\alpha}(\Omega)},
\|v\|_{C^{1,\alpha}(\Omega)}\leq C
\end{equation}
for some constant $C>0$, and the following asymptotic behavior in
far fields
\begin{eqnarray}
\rho\rightarrow \rho_0>0,  \,\, u\rightarrow u_0(x_2)>0, \,\,
v\rightarrow 0\,\, \text{as}\,\, x_1\rightarrow -\infty,\label{Easymptoticupstream1}\\
\nabla \rho \rightarrow 0,\,\,\nabla u\rightarrow (0, u_0'(x_2)),
\,\, \nabla v\rightarrow 0\,\, \text{as}\,\, x_1\rightarrow
-\infty,\label{Easymptoticupstream2}
\end{eqnarray}
uniformly for $x_2\in K_1\Subset(0,1)$, and
\begin{eqnarray}
\rho\rightarrow \rho_1>0,  \,\, u\rightarrow u_1(x_2)>0, \,\,
v\rightarrow 0\,\, \text{as}\,\, x_1\rightarrow
+\infty,\label{Easymptoticdownstream1}\\ \nabla \rho \rightarrow
0,\,\,\nabla u\rightarrow (0, u_1'(x_2)), \,\, \nabla v\rightarrow
0\,\, \text{as}\,\, x_1\rightarrow
+\infty,\label{Easymptoticdownstream2}
\end{eqnarray}
uniformly for $x_2\in K_2\Subset (a, b)$, where $\rho_0$ and
$\rho_1$ are both positive constants, and $\rho_0$, $\rho_1$, $u_0$,
and $u_1$ can be determined by $m$, $B(x_2)$ and $b-a$ uniquely;
\item (Uniqueness) the Euler flow which satisfies
(\ref{EEulercontinuityeq})-(\ref{EEulermomentumeq2}), boundary
condition (\ref{Enoflowbc}), asymptotic condition
(\ref{EasymtoticBernoulli}), mass flux condition (\ref{Emassflux}),
(\ref{Esubsonic}), and asymptotic behavior
(\ref{Easymptoticupstream1})-(\ref{Easymptoticupstream2}) is unique;
\item (Critical mass flux) If, besides
(\ref{EassumptiononBernoulli}), $B$ also satisfies
\begin{equation}\label{Ecricondition}
B'(0)=B'(1)=0,
\end{equation}
then $\hat{m}$ is the upper critical mass flux for the existence of
subsonic flow in the following sense: either
\begin{equation}
\sup_{\bar{\Omega}}(u^2+v^2-c^2(\rho))\rightarrow
0\,\,\text{as}\,\, m\rightarrow \hat{m},
\end{equation}
or there is no $\sigma>0$ such that for all $m\in
(\hat{m},\hat{m}+\sigma)$, there are Euler flows satisfying
(\ref{EEulercontinuityeq})-(\ref{EEulermomentumeq2}), boundary
condition (\ref{Enoflowbc}), asymptotic condition
(\ref{EasymtoticBernoulli}), mass flux condition
(\ref{Emassflux}), (\ref{Esubsonic}), and asymptotic behavior
(\ref{Easymptoticupstream1})-(\ref{Easymptoticupstream2}) and
\begin{equation}
\sup_{m\in(\hat{m},\hat{m}+\sigma)}\sup_{\bar{\Omega}}(c^2(\rho)-(u^2+v^2))>0.
\end{equation}
\end{enumerate}
\end{thm}

There are a few remarks in order:
\begin{rmk}
{\rm  Here we obtained only the existence of the Euler flows in
the nozzle, and the uniqueness in a special class of flows, but
not the uniqueness for general Euler flows. For the issue on the
uniqueness for steady incompressible Euler flows in a bounded
domain, please refer to \cite{Troshkin}.}
\end{rmk}

\begin{rmk}
{\rm It can be shown by modifying the analysis in this paper
slightly without further difficulties that there exists a subsonic
full compressible Euler flow in the nozzle, if the entropy is
prescribed in the upstream.}
\end{rmk}

\begin{rmk}
{\rm The subsonic Euler flows in half plane was studied in
\cite{ChenJun} recently. Although stream function formulation is
also introduced in \cite{ChenJun}, however, the far fields
conditions are different from ours. Furthermore, we obtain
critical upper bound of mass flux for existence of subsonic flows
in nozzles.}
\end{rmk}

The rest of the paper is arranged as follows: in Section
\ref{SEexistence}, we reformulate the problem by deriving the
governing equation and boundary conditions for Euler flows in
terms of a stream function, provided that the Euler flow has
simple topological structure and satisfies the asymptotic behavior
(\ref{Easymptoticupstream1})-(\ref{Easymptoticupstream2}). In
Section \ref{SEexistence}, existence of solutions to a modified
elliptic problem is established. Subsequently, in Section
\ref{SEasymptotic}, we will study asymptotic behavior of solutions
in a larger class and show uniqueness of the solution to the
boundary value problem. The existence of boundary value problem
for the stream functions will be a direct consequence of these
asymptotic behavior and uniqueness.  In Section \ref{SErefined},
some refined estimates for the stream function will be derived.
Combining these estimates with the asymptotic behavior obtained in
Section \ref{SEasymptotic} will yield the existence of Euler flows
which satisfy all properties in Theorem \ref{EThexistence}.
Finally, in Section \ref{SEcritical}, we will show the existence
of the critical mass flux.

\section{Stream-Function Formulation of the Problem}\label{SEreformulation}
We start with some basic structures of the steady Euler system. The
steady Euler system
(\ref{EEulercontinuityeq})-(\ref{EEulermomentumeq2}) can be written
in the following form,
\begin{equation*}
AU_{x_1}+BU_{x_2}=0,
\end{equation*}
where
\begin{equation*}
A=
\begin{pmatrix}
\frac{uc^2(\rho)}{\rho} & c^2(\rho) & 0\\
c^2(\rho) & \rho u  & 0\\
0 & 0 & \rho u
\end{pmatrix}
, \,\, B=
\begin{pmatrix}
\frac{vc^2(\rho)}{\rho} & 0 & c^2(\rho)\\
0 & \rho v  & 0\\
c^2(\rho) & 0 & \rho v
\end{pmatrix}
, \,\,U=
\begin{pmatrix}
\rho\\
u\\
v
\end{pmatrix}
.
\end{equation*}
Let $\lambda$ be the solution of
\begin{equation}\label{Eeigenvalueeq}
\det (\lambda A-B)=0.
\end{equation}
It follows from straightforward computations that
(\ref{Eeigenvalueeq}) has three eigenvalues
\begin{equation*}
\lambda_1=\frac{v}{u},\,\lambda_{\pm}=\frac{uv\pm
c(\rho)\sqrt{u^2+v^2-c^2(\rho)}}{u^2-c^2}.
\end{equation*}
Therefore, at the points where $u^2+v^2-c^2(\rho)>0$, i.e., the
flow is supersonic, (\ref{Eeigenvalueeq}) has 3 real eigenvalues,
the Euler system is hyperbolic. When $u^2+v^2-c^2(\rho)<0$, i,e.,
the flow is subsonic, (\ref{Eeigenvalueeq}) has a real eigenvalue
and two complex eigenvalues, the Euler system is a
hyperbolic-elliptic coupled system. Therefore, even for globally
subsonic flows, one has to resolve a hyperbolic mode. Moreover,
for flows in infinitely long nozzles with both ends at infinity,
it seems difficult to get uniform estimates for hyperbolic mode.

To overcome the difficulties mentioned above, we introduce the
stream functions for the 2-D steady compressible Euler flows, and
derive an equivalent formulation for Euler flows in terms of the
stream functions when the flow satisfies certain asymptotic
behavior.

It follows from (\ref{Enonconservemeq1}) and
(\ref{Enonconservemeq2}) that
\begin{equation}\label{Emixedmomentumeq}
\partial_{x_2}(uu_{x_1}+vu_{x_2}+h(\rho)_{x_1})
-\partial_{x_1}(uv_{x_1}+vv_{x_2}+h(\rho)_{x_2})=0,
\end{equation}
which yields that
\begin{equation}\label{Evorticityteq1}
(u,v)\cdot\nabla \omega+\omega\Div(u,v)=0,
\end{equation}
where $\omega=v_{x_1}-u_{x_2}$ is the vorticity of the flow. By the
continuity equation (\ref{EEulercontinuityeq}), one has
\begin{equation*}
(u,v)\cdot\nabla \omega+\omega\Div(u,v)=(u,v) \cdot(\nabla
\omega-\frac{\omega\nabla\rho}{\rho})=\rho(u,v)\cdot\nabla\left(\frac{\omega}{\rho}\right).
\end{equation*}
Therefore, away from vacuum, (\ref{Evorticityteq1}) is equivalent to
\begin{equation}\label{Evorticityeq}
(u,v)\cdot\nabla\left(\frac{\omega}{\rho}\right)=0.
\end{equation}
We have the following proposition.
\begin{prop}\label{Elemmaequivalent}
For a smooth flow away from vacuum in the nozzle $\Omega$
satisfying (\ref{Enontrivialbdy}) and (\ref{Eboundary1}), the
system consisting of (\ref{EEulercontinuityeq}),
(\ref{EstreamBernoulli}) and (\ref{Evorticityeq}) is equivalent to
the original Euler equations
(\ref{EEulercontinuityeq})-(\ref{EEulermomentumeq2}), if the given
flow satisfies no flow boundary condition,
\begin{equation}\label{Ephorizontalvelocity}
u>0\,\, \text{in}\,\,\Omega,
\end{equation}
and the following asymptotic behavior
\begin{equation}\label{EcequivstreamEuler}
u,\,\,\rho\,\, \text{and}\,\, v_{x_2}\,\, \text{are bounded,
 while}\,\, v,\,\,v_{x_1}\,\,\text{and}\,\,\rho_{x_2}\rightarrow
0,\,\,\text{as}\,\,x_1\rightarrow -\infty.
\end{equation}
\end{prop}
\begin{pf}
From previous analysis, it is easy to see that smooth solutions to
the Euler equations
(\ref{EEulercontinuityeq})-(\ref{EEulermomentumeq2}) satisfy
(\ref{EEulercontinuityeq}), (\ref{EstreamBernoulli}) and
(\ref{Evorticityeq}). On the other hand, it follows from
(\ref{EEulercontinuityeq}), (\ref{Evorticityeq}) and the above
derivation that (\ref{Emixedmomentumeq}) holds. Therefore, there
exists a function $\Phi$ such that
\begin{equation*}
\Phi_{x_1}=uu_{x_1}+vu_{x_2}+h(\rho)_{x_1}, \,\,
\Phi_{x_2}=uv_{x_1}+vv_{x_2}+h(\rho)_{x_2}.
\end{equation*}
So, (\ref{EstreamBernoulli}) is equivalent to
\begin{equation}\label{EtransportPotential}
(u,v)\cdot\nabla\Phi=0.
\end{equation}

Due to the no flow boundary condition (\ref{Enoflowbc}), $\Phi$ is a
constant along each component of the nozzle boundary. If, in
addition,
\begin{equation}\label{EupstramPotentialderivative}
\Phi_{x_2}\rightarrow 0 \,\, \text{as}\,\, x_1\rightarrow -\infty,
\end{equation}
then $\Phi\rightarrow C$ as $x_1\rightarrow -\infty$. On the other
hand, it follows from (\ref{Ephorizontalvelocity}) that through
each point in $\Omega$, there is one and only one streamline
satisfying
\begin{equation*}
\left\{
\begin{array}{ll}
\frac{dx_1}{ds}=u(x_1(s),x_2(s)),\\
\frac{dx_2}{ds}=v(x_1(s),x_2(s)),
\end{array}
\right.
\end{equation*}
which can be defined globally in the nozzle (i.e., from the entry
to the exit). Furthermore, it follows from
(\ref{EEulercontinuityeq}) that any streamline through some point
in $\Omega$ can not touch the nozzle wall. Suppose not, let the
streamline through $(x_1^0,x_2^0)$ pass through $(x_1, f_1(x_1))$.
Due to (\ref{EEulercontinuityeq}) and no flow boundary condition,
one has
\begin{equation*}
\int_{f_1(x_1^0)}^{x_2^0}(\rho u)(x_1^0,x_2)dx_2=0.
\end{equation*}
This contradicts (\ref{Ephorizontalvelocity}).

Thus, one can always solve (\ref{EtransportPotential}) in the whole
domain $\Omega$, which yields
\begin{equation*}
\Phi\equiv C\,\,\text{in}\,\,\Omega,
\end{equation*}
if (\ref{EupstramPotentialderivative}) holds. Therefore,
$\Phi_{x_1}=\Phi_{x_2}\equiv 0$ in $\Omega$, i.e.,
(\ref{Enonconservemeq1}) and (\ref{Enonconservemeq2}) hold
globally in the nozzle. Thus, both (\ref{EEulermomentumeq1}) and
(\ref{EEulermomentumeq2}) are true. It is obvious that
(\ref{EupstramPotentialderivative}) holds if
(\ref{EcequivstreamEuler}) is valid.
\end{pf}

It suffices to prove the existence of solutions to the system
(\ref{EEulercontinuityeq}), (\ref{EstreamBernoulli}) and
(\ref{Evorticityeq}) satisfying (\ref{Ephorizontalvelocity}) and
(\ref{EcequivstreamEuler}).

However, system (\ref{EEulercontinuityeq}),
(\ref{EstreamBernoulli}) and (\ref{Evorticityeq}) is not easy to
study directly either, since for infinitely long nozzles with both
ends at infinity, it seems difficult to estimate the solutions for
transport equations (\ref{EstreamBernoulli}) and
(\ref{Evorticityeq}). Instead, we will use an equivalent
formulation for (\ref{EEulercontinuityeq}),
(\ref{EstreamBernoulli}) and (\ref{Evorticityeq}).

It follows from (\ref{EEulercontinuityeq}) that there exists a
stream function $\psi$ such that
\begin{equation*}
\psi_{x_1}=-\rho v, \,\, \psi_{x_2}=\rho u.
\end{equation*}
Thus for the flows away from the vacuum, (\ref{Evorticityeq}) is
equivalent to
\begin{equation}\label{Estreamvorticitytraneq}
\nabla^{\bot}\psi\cdot \nabla\left(\frac{\omega}{\rho}\right)=0,
\end{equation}
where $\nabla^{\bot}=(-\partial_{x_2},\partial_{x_1})$. Note that
(\ref{Estreamvorticitytraneq}) means that $\frac{\omega}{\rho}$
and $\psi$ are functionally dependent, therefore, one may regard
$\frac{\omega}{\rho}$ as a function of $\psi$. Set
\begin{equation}\label{Evorticitystreamrelation}
\frac{\omega}{\rho}=W(\psi).
\end{equation}
 Similarly, (\ref{EstreamBernoulli}) is equivalent to
\begin{equation*}
\nabla^{\bot}\psi\cdot\nabla(h(\rho)+\frac{1}{2}(u^2+v^2))=0,
\end{equation*}
therefore, $h(\rho)+\frac{1}{2}(u^2+v^2)$ is also a function of
$\psi$. We define this function by
\begin{equation}\label{EBernoullistreamrelation}
h(\rho)+\frac{1}{2}(u^2+v^2)=\mathcal{B}(\psi).
\end{equation}
Furthermore, it follows from (\ref{Enoflowbc}) that the nozzle walls
are streamlines, so $\psi$ is constant on each nozzle wall. Due to
(\ref{Emassflux}), one can assume that
\begin{equation}\label{Estreambc}
\psi=0\,\,\text{on} \,\,S_1,\,\,\text{and}\,\, \psi=m\,\,\text{on}
\,\,S_2.
\end{equation}

In order to get an explicit form of $\mathcal{B}$, first we study
the density-speed relation via Bernoulli's law
(\ref{EBernoullistreamrelation}) carefully.

Note that $p'(\rho)>0$ for $\rho>0$ and $p''(\rho)\geq 0$,
therefore $h'(\rho)=p'(\rho)/\rho>0$ for $\rho>0$ and for some
fixed $\bar{\rho}>0$,
\begin{equation*}
h(\rho)=h(\bar{\rho})+\int_{\bar{\rho}}^{\rho}\frac{p'(s)}{s}ds\geq
h(\bar{\rho})+\int_{\bar{\rho}}^{\rho}\frac{p'(\bar{\rho})}{s}ds\,\,
\text{for} \,\,\rho>\bar{\rho}.
\end{equation*}
This yields that $h(\rho)\rightarrow \infty$ as $\rho\rightarrow
\infty$. On the other hand, since $\inf_{\rho>0}h(\rho)=B_0$,
$h(\rho)\rightarrow B_0$ as $\rho\rightarrow 0$. Thus for any
$s>B_0$, there exists a unique $\bar{\varrho}=\bar{\varrho}(s)>0$
such that
\begin{equation*}
h(\bar{\varrho}(s))=s.
\end{equation*}
Moreover, for the state with given Bernoulli's constant $s$, the
density and speed satisfy the relation,
\begin{equation*}
h(\rho)+\frac{q^2}{2}=s.
\end{equation*}
Therefore, the speed $q$ satisfies
\begin{equation*}
q=\sqrt{2(s-h(\rho))}.
\end{equation*}
Hence, for fixed $s$, $q$ is a strictly decreasing function of
$\rho$ on $[0,\bar{\varrho}(s)]$. By the definition of
$\bar{\varrho}(s)$, one has
$q(\bar{\varrho}(s))=0<c(\bar{\varrho}(s))$. Now we claim that
$q(0)>c(0)$. Indeed, one can prove this claim in two cases. First,
if $c(0)>0$, then
\begin{equation*}
h(\rho)=h(\bar{\rho})+\int_{\bar{\rho}}^{\rho}\frac{p'(s)}{s}ds\leq
h(\bar{\rho})+\int_{\bar{\rho}}^{\rho}\frac{c^2(0)}{s}ds
\rightarrow -\infty\,\,\text{as}\,\,\rho\rightarrow 0.
\end{equation*}
Therefore, $q(0)\rightarrow \infty$. Thus $q(0)>c(0)$. Second, if
$c(0)=0$, $q(\rho)\rightarrow \sqrt{2(s-B_0)}$ as $\rho\rightarrow
0$, therefore, $q(0)>0=c(0)$. This completes the proof of the
claim. Since $c^2(\rho)=p'(\rho)$ is an increasing function of
$\rho$, there exists a unique $\varrho(s)\in [0,\bar{\varrho}(s)]$
such that
\begin{equation*}
c^2(\varrho(s))=q^2(\varrho(s)).
\end{equation*}
In summary, for any given $s>B_0$, there exist
$\bar{\varrho}=\bar{\varrho}(s)$, $\varrho=\varrho(s)$ and
$\Gamma=\Gamma(s)$ such that
\begin{equation}\label{Edefmaxcrdensity}
h(\bar{\varrho}(s))=s,\,\,
h(\varrho(s))+\frac{\Gamma^2(s)}{2}=s,\,\,\text{and}\,\,
c^2(\varrho(s))=\Gamma^2(s),
\end{equation}
where $\bar{\varrho}(s)$, $\varrho(s)$, and $\Gamma(s)$ are the
maximum density, the critical density, and the critical speed,
respectively for the states with given Bernoulli's constant $s$.
Set
\begin{equation}\label{Edefcirticalmassfluxfunction}
\Sigma(s)=\varrho(s)\sqrt{2(s-h(\varrho(s)))}.
\end{equation}
Then direct calculations show that
\begin{equation*}
\frac{d \bar{\varrho}}{d
s}=\frac{\bar{\varrho}}{p'(\bar{\varrho})},\,\,\frac{d
\varrho}{ds}=\frac{1}{\frac{p'(\varrho)}{\varrho}+\frac{p''(\varrho)}{2}},
\end{equation*}
and
\begin{equation*}
\frac{d \Sigma}{ds}=
\frac{\sqrt{2(s-h(\varrho(s)))}}{\frac{p'(\varrho)}{\varrho}+\frac{p''(\varrho)}{2}}
+\varrho \frac{1-\frac{2p'(\varrho)}{2 p'(\varrho)+\varrho
p''(\varrho)}} {\sqrt{2(s-h(\varrho(s)))}}.
\end{equation*}
Thus
\begin{equation*}
\frac{d \bar{\varrho}}{d s}>0,\,\,\frac{d
\varrho}{ds}>0,\,\,\text{and}\,\, \frac{d\Sigma}{ds}>0.
\end{equation*}
Obviously, $\varrho(s)<\bar{\varrho}(s)$, if $s>B_0$. By the
continuity and monotonicity of $\varrho(s)$ and $\bar{\varrho}(s)$,
there exists a unique $\underline{\delta}>0$ such that
\begin{equation}\label{EconditionforB}
\varrho(\underline{B}+\underline{\delta})=\bar{\varrho}(\underline{B}).
\end{equation}
Moreover, it follows from (\ref{Edefmaxcrdensity}) that there exists
a uniform constant $C>0$ such that
\begin{equation}\label{EuniformestimateBernoulli}
\left\{
\begin{array}{ll}
 C^{-1}\leq
\varrho(\underline{B})<\bar{\varrho}(\underline{B})=\varrho(\underline{B}+\underline{\delta})
<\bar{\varrho}(\underline{B}+\underline{\delta})\leq
C,\\
C^{-1}\leq \varrho'(s)\leq C,\,\,C^{-1}\leq \bar{\varrho}'(s)\leq C,
\,\,\text{if}\,\,s\in(\underline{B},\underline{B}+\underline{\delta}),\\
C^{-1}\leq h'(\rho)\leq
C,\,\,\text{if}\,\,\rho\in(\varrho(\underline{B}),\bar{\varrho}(\underline{B}+\underline{\delta})),\\
C^{-1}\leq \Sigma(s)\leq C,\,\,\text{if}\,\, s\in
(\underline{B},\underline{B}+\underline{\delta}).
\end{array}
\right.
\end{equation}
Later on,  $C$ will  to denote generic constants which depend only
on $\underline{B}$ and $\underline{\delta}$, and thus essentially on
$\underline{B}$.

In order to study the relationship between density and mass flux
with given Bernoulli's constant, let us investigate the function
defined by
\begin{equation*}
I(\rho)=2\rho^2(s-h(\rho)).
\end{equation*}
Direct calculations show
\begin{equation*}
\frac{d
I}{d\rho}=4\rho(s-h(\rho)-p'(\rho)/2)=2\rho(q^2(\rho)-c^2(\rho)).
\end{equation*}
Therefore, for $\rho\in (0,\varrho(s))$,  $\frac{d I}{d\rho}>0$; and
$\frac{d I}{d\rho}<0$ for $\rho\in(\varrho(s),\bar{\varrho}(s))$.
Moreover, $I(0)=I(\bar{\varrho}(s))=0$. Thus $I(\rho)>0$ if $\rho\in
(0,\bar{\varrho}(s))$ and $I$ achieves its maximum at
$\rho=\varrho(s)$. So, for fixed $s$, the relation
\begin{equation}
h(\rho)+\frac{\mathcal{M}}{2\rho^2}=s
\end{equation}
defines a function $\mathcal{M}=I(\rho)$ which attains its maximum
$\mathcal{M}=\Sigma^2(s)$ at $\rho=\varrho(s)$. Thus $\rho$ is a
two-valued function of $\mathcal{M}$ for $\mathcal{M}\in
[0,\Sigma^2(s))$. Denote the subsonic branch by
\begin{equation*}
\rho=J(\mathcal{M})\,\,\text{for}\,\,\mathcal{M}\in(0,\Sigma^2(s)),
\end{equation*}
which satisfies $J(\mathcal{M})>\varrho(s)$. When $s$ varies,  this
branch will be denoted by
\begin{equation}\label{Edefsubsonicbranch}
\rho=J(\mathcal{M},s)\,\,\text{for}\,\,(\mathcal{M},s)\in
\{(\mathcal{M},s)|\mathcal{M}\in(0,\Sigma^2(s)),s>B_0\}.
\end{equation}

To determine the explicit form of $W$ and $\mathcal{B}$, one may
study  $W$ and $\mathcal{B}$ in the far fields of the nozzle where
the flow may have certain simple asymptotic structure. Indeed, for
flows satisfying the asymptotic behavior
(\ref{Easymptoticupstream1})-(\ref{Easymptoticdownstream2}), one can
determine $\rho_0$, $\rho_1$, $u_0(x_2)$ and $u_1(x_2)$ first.
Suppose that the flow satisfies (\ref{Easymptoticupstream1}). Then
\begin{equation}\label{EupstreamBernoulli}
h(\rho_0)+\frac{u_0^2(x_2)}{2}=B(x_2),\,\, u_0(x_2)> 0,
\end{equation}
and
\begin{equation}\label{Eupstreammassflux}
\int_0^1\rho_0u_0(x_2)dx_2=m
\end{equation}
hold, which shows that
\begin{equation}\label{Eupstreamdensityspeed}
u_0(x_2)=\sqrt{2(B(x_2)-h(\rho_0))},
\end{equation}
and
\begin{equation}\label{Easymtoticmassflux}
m=\int_0^1\rho_0\sqrt{2(B(x_2)-h(\rho_0))}dx_2.
\end{equation}
If $B(x_2)$ satisfies
\begin{equation}\label{EconditionBagain}
\inf_{x_2\in[0,1]}B(x_2)=\underline{B},\,\,\|B'(x_2)\|_{C^{0,1}([0,1])}\leq
\delta,
\end{equation}
then
\begin{equation}
\bar{B}=\sup_{x_2\in[0,1]}B(x_2)\leq \underline{B}+\delta.
\end{equation}
Let $\delta\leq \underline{\delta}/2$. Then it follows from
(\ref{EconditionforB}) that $\varrho(B(x_2))\leq
\varrho(\bar{B})<\bar{\varrho}(\underline{B})$. To obtain a global
subsonic flow in the nozzle, it is necessary to show that  for given
$B_2(x_2)$ and $m$, (\ref{Easymtoticmassflux}) has a solution
satisfying $\rho_0\in
(\varrho(\bar{B}),\bar{\varrho}(\underline{B}))$. Direct
calculations yield that
\begin{equation*}
\frac{d}{d\rho_0}\int_0^1\rho_0\sqrt{2(B(x_2)-h(\rho_0))}dx_2<0,
\,\, \text{for}\,\,
\rho_0\in(\varrho(\bar{B}),\bar{\varrho}(\underline{B})).
\end{equation*}
It follows from (\ref{EuniformestimateBernoulli}) and
(\ref{EconditionBagain}) that
\begin{equation*}
\int_{0}^1
\bar{\varrho}(\underline{B})\sqrt{2(B(x_2)-h(\bar{\varrho}(\underline{B})))}dx_2
= \int_0^1
\bar{\varrho}(\underline{B})\sqrt{2(B(x_2)-\underline{B})}dx_2
\leq C \delta^{1/2}.
\end{equation*}
In addition,
\begin{eqnarray*}
&&\int_{0}^1
\varrho(\bar{B})\sqrt{2(B(x_2)-h(\varrho(\bar{B})))}dx_2\\
&\geq& \int_0^1
\varrho(\bar{B})\sqrt{2(\underline{B}-h(\varrho(\bar{B})))}dx_2\\
&=& \int_0^1
\varrho(\bar{B})\sqrt{2(h(\bar{\varrho}(\underline{B}))-h(\varrho(\bar{B})))}dx_2\\
&=& \int_0^1
\varrho(\bar{B})\sqrt{2(h(\varrho(\underline{B}+\underline{\delta}))-h(\varrho(\bar{B})))}dx_2\\
&\geq & \int_0^1
\varrho(\bar{B})\sqrt{2\left(h(\varrho(\underline{B}+
\underline{\delta}))-h(\varrho(\underline{B}+\underline{\delta}/2))\right)}dx_2\\
&\geq& C^{-1}\underline{\delta}^{1/2}.
\end{eqnarray*}
Therefore, for any $\gamma\in(0,1/3)$, there exists
$\tilde{\delta}_0\in (0,  \underline{\delta}/2)$ such that
(\ref{Easymtoticmassflux}) admits a unique solution
$\rho_0\in(\varrho(\bar{B}),\bar{\varrho}(\underline{B}))$, provided
that $0\leq \delta\leq \tilde{\delta}_0$ and $m\in (\delta^{\gamma},
m_1)$, where $m_1$ satisfies $m_1\geq
C^{-1}\underline{\delta}^{1/2}\geq 2\tilde{\delta}_0^{\gamma/2}$.
Later on, for definiteness, we will choose $\gamma=1/4$. However,
all results hold for $\gamma\in (0,1/3)$.

By virtue of (\ref{Easymtoticmassflux}), one has
\begin{eqnarray*}
m&=&\int_0^1\rho_0\sqrt{2(B(x_2)-h(\rho_0))}dx_2\\
&=&
\int_0^1\rho_0\sqrt{2(B(x_2)-\underline{B}+\underline{B}-h(\rho_0))}dx_2\\
&\leq&
C\int_0^1\sqrt{2(\delta+h(\bar{\varrho}(\underline{B}))-h(\rho_0))}dx_2.
\end{eqnarray*}
Thus
\begin{equation*}
\delta+h(\bar{\varrho}(\underline{B}))-h(\rho_0)\geq  C^{-1} m^2\geq
C^{-1} \delta^{2\gamma}.
\end{equation*}
Note that $\gamma<1/3$, therefore, there exists
$\tilde{\tilde{\delta}}_0\in (0, \tilde{\delta}_0)$ such that if
$0<\delta\leq \tilde{\tilde{\delta}}_0$, then
\begin{equation*}
h(\bar{\varrho}(\underline{B}))-h(\rho_0)\geq
C^{-1}\delta^{2\gamma}.
\end{equation*}
Consequently, if $\|B'(x_2)\|_{C^{0,1}([0,1])}= \delta\leq
\tilde{\tilde{\delta}}_0$, by virtue of
(\ref{EuniformestimateBernoulli}), there is a positive constant
$C$ such that
\begin{equation}\label{Efarestimatebystream}
\left\{
\begin{array}{ll}
C^{-1}\delta^{2\gamma}\leq \bar{\varrho}(\underline{B})-\rho_0\leq C,\\
C^{-1}\delta^{\gamma}\leq u_0\leq C,\\
|u_0'(x_2)|=\left|\frac{B'(x_2)}{\sqrt{2(B(x_2)-h(\rho_0))}}\right|\leq
C \delta^{1-\gamma},\\[1mm]
[u_0'(x_2)]_{C^{0,1}([0,1])}\leq
C(\delta^{1-\gamma}+\delta^{2-3\gamma}).
\end{array}
\right.
\end{equation}

To determine the states in the downstream, we parametrize the
streamlines in the downstream by their positions in the upstream.
Due to (\ref{Easymptoticupstream1}), (\ref{Easymptoticdownstream1}),
and (\ref{Ephorizontalvelocity}), we can define
\begin{equation}\label{Eflowmap}
y=y(s)\,\,\text{for}\,\,s\in[0,1]
\end{equation}
such that
\begin{eqnarray}
&&h(\rho_0)+\frac{u_0^2(s)}{2}=h(\rho_1)+\frac{u_1^2(y(s))}{2},\,\,
u_1(y(s))>0,\label{EfarBernoulliconservation}\\
&&\int_0^s\rho_0u_0(t)dt=\int_a^{y(s)}\rho_1u_1(t)dt,\label{Efarmassconservation}\\
&&y(0)=a,\,\, y(1)=b.\label{EflowmapBC}
\end{eqnarray}
The meaning of $y(s)$ is that the streamline which starts from
$(-\infty,s)$ will flow to $(\infty,y(s))$. The map
(\ref{Eflowmap}) is well-defined since
(\ref{Ephorizontalvelocity}) ensures a simple topological
structure of streamlines. It follows from
(\ref{Efarmassconservation}) that
\begin{equation}\label{Efarmasscondiff}
\rho_0 u_0(s)=\rho_1 u_1(y(s))y'(s).
\end{equation}
Hence,
\begin{equation}\label{EfarparaODE}
\left\{
\begin{array}{ll}
\frac{dy}{ds}=\frac{\rho_0u_0(s)}{\rho_1
\sqrt{2(h(\rho_0)+\frac{u_0^2(s)}{2}-h(\rho_1))}},\\
y(0)=a,
\end{array}
\right.
\end{equation}
where the parameter $\rho_1$ satisfies
\begin{equation}\label{Efarparaheight}
\int_0^1\frac {\rho_0u_0(s)}{\rho_1
\sqrt{2(h(\rho_0)+\frac{u_0^2(s)}{2}-h(\rho_1))}}ds=b-a.
\end{equation}
It remains to show that there exists a $\rho_1\in
(\varrho(\bar{B}),\bar{\varrho}(\underline{B}))$ satisfying
(\ref{Efarparaheight}). By direct calculations, one has
\begin{equation*}
\frac{d}{d\rho_1}\int_0^1\frac {\rho_0u_0(s)}{\rho_1
\sqrt{2(h(\rho_0)+\frac{u_0^2(s)}{2}-h(\rho_1))}}ds> 0,\,\,
\text{for}\,\,
\rho_1\in(\varrho(\bar{B}),\bar{\varrho}(\underline{B})).
\end{equation*}
First, there exists $\bar{\delta}_0\in (0,
\tilde{\tilde{\delta}}_0)$ such that
\begin{eqnarray*}
&&\int_0^1\frac {\rho_0u_0(s)}{\bar{\varrho}(\underline{B}-\delta)
\sqrt{2(h(\rho_0)+\frac{u_0^2(s)}{2}-h(\bar{\varrho}(\underline{B}-\delta)))}}ds\\
& = & \int_0^1\frac
{\rho_0u_0(s)}{\bar{\varrho}(\underline{B}-\delta)
\sqrt{2(B(s)-\underline{B}+h(\bar{\varrho}(\underline{B}))
-h(\bar{\varrho}(\underline{B}-\delta)))}}ds\\
&\geq& C\delta^{(2\gamma-1)/2}>b-a,
\end{eqnarray*}
provided $\delta\leq \bar{\delta}_0$. On the other hand,
\begin{eqnarray*}
&&\int_0^1\frac {\rho_0u_0(s)}{\varrho(\bar{B})
\sqrt{2(h(\rho_0)+\frac{u_0^2(s)}{2}-h(\varrho(\bar{B})))}}ds\\
&\leq& \frac{m}{\varrho(\bar{B})
\sqrt{2(\underline{B}-h(\varrho(\bar{B})))}}\\
&=& \frac{m}{\varrho(\bar{B})
\sqrt{2(h(\bar{\varrho}(\underline{B}))-h(\varrho(\bar{B})))}}\\
& = & \frac{m}{\varrho(\bar{B})
\sqrt{2\left(h(\varrho(\underline{B}+\underline{\delta}))-h(\varrho(\bar{B}))\right)}}\\
&\leq & \frac{m}{\varrho(\bar{B})
\sqrt{2\left(h(\varrho(\underline{B}+\underline{\delta}))-h(\varrho(\underline{B}+\underline{\delta}/2))\right)}}\\
 &\leq & \frac{m}{C^{-1}\underline{\delta}^{1/2}}.
\end{eqnarray*}
So there exists a unique $\rho_1\in
(\varrho(\underline{B}),\bar{\varrho}(\underline{B}))$ such that
(\ref{Efarparaheight}) holds, if $0\leq \delta\leq\bar{\delta}_0$
and $m\in (\delta^{\gamma},m_2)$ for some $m_2\geq \min\{m_1,
C^{-1}(b-a)\underline{\delta}^{1/2}\}$. Furthermore, one can choose
$\bar{\delta}_0$ smaller if necessary such that $m_2\geq
2\bar{\delta}_0^{\gamma/2}$. As soon as $\rho_1$ is determined,
$y(s)$ and $u_1$ can be obtained from (\ref{EfarparaODE}) and
(\ref{EfarBernoulliconservation}).

Let us summarize the above calculations in the following
proposition:
\begin{prop}
Let $\underline{B}>B_0$. There exists $\bar{\delta}_0>0$ such that
for any $B\in C^{1,1}([0,1])$ satisfying (\ref{EconditionBagain})
with $\delta\leq \bar{\delta}_0$, there exists $\bar{m}\geq
2\bar{\delta}_0^{1/8}$ such that
\begin{enumerate}
\item
 there exist solutions $(\rho_0,u_0)$
to (\ref{EupstreamBernoulli})-(\ref{Eupstreammassflux})  and
$(\rho_1,u_1)$ solving
(\ref{EfarBernoulliconservation})-(\ref{EflowmapBC}) if
$m\in(\delta^{1/4},\bar{m})$; \item $\rho_0$, $\rho_1\in
(\varrho(\bar{B}),\bar{\varrho}(\underline{B}))$; \item $(\rho_0,
u_0)$ satisfies (\ref{Efarestimatebystream});
\item
 either
$\rho_0\rightarrow \varrho(\bar{B})$ or $\rho_1\rightarrow
\varrho(\bar{B})$ as $m\rightarrow \bar{m}$;
\end{enumerate}
 where
$\bar{B}=\max_{x_2\in[0,1]}B(x_2)$.
\end{prop}
\begin{pf}
We need only to verify the last statement.

First, if $m\in (\delta^{1/4},m_2)$, both $\rho_0$ and $\rho_1$
belong to $(\varrho(\bar{B}),\bar{\varrho}(\underline{B}))$. For a
given $B(x_2)$, as $m$ increases, $\rho_0$ decreases. If
$m\rightarrow\tilde{m}=\int_0^1\varrho(\bar{B})\sqrt{2(B(x_2)-h(\varrho(\bar{B})))}dx_2$,
then $\rho_0\rightarrow\varrho(\bar{B})$. Therefore there exists an
upper bound for $m$ to ensure the existence of $\rho_0$, $\rho_1\in
(\varrho(\bar{B}),\bar{\varrho}(\underline{B}))$. Define
\begin{equation}\label{Edefmaxflux}
\bar{m}=\sup\{s|m\in(\delta^{\gamma},s),\text{ there exist  }\rho_0,
\rho_1\in (\varrho(\bar{B}),\bar{\varrho}(\underline{B}))\}
\end{equation}
Obviously, $\bar{m}\in [m_2,\tilde{m}]$. Note that $\rho_0$ and
$\rho_1$ are uniformly away from $\bar{\varrho}(\underline{B})$. If
neither $\rho_0$ nor $\rho_1$ approaches to $\varrho(\bar{B})$ as
$m\rightarrow\bar{m}$, then there always exist $\rho_0$, $\rho_1\in
(\varrho(\bar{B}),\bar{\varrho}(\underline{B}))$ for $m\in
(\delta^{\gamma},\bar{m}+\epsilon)$ with some small $\epsilon>0$.
This contradicts with the definition of $\bar{m}$. So the proof of
the Proposition is finished.
\end{pf}

We now can determine $W$ and $\mathcal{B}$ in the upstream.
Suppose that the flow satisfies the asymptotic behavior
(\ref{Easymptoticupstream1}). Then in the upstream, a stream
function can be chosen so that
\begin{equation}
\psi=\int_0^{X_2}\rho_0u_0(s)ds,
\end{equation}
and $0\leq \psi\leq m$. Since $\rho_0u_0(s)>0$ for $s\in[0,1]$,
$\psi$ is an increasing function of $X_2$. Thus one can represent
$X_2$ as a function of $\psi$,
\begin{equation*}
X_2=\kappa (\psi),\,\,\,\, 0\leq \psi\leq m.
\end{equation*}
Define
\begin{equation}\label{Eupstreamviastream}
f(\psi)=u_0'(\kappa(\psi)),\,\,\text{and}\,\,
F(\psi)=u_0(\kappa(\psi)).
\end{equation}
Then $f$ and $F$ are well-defined on $[0,m]$.

\begin{center}
     \includegraphics[width=9.1cm,height=3.1cm]{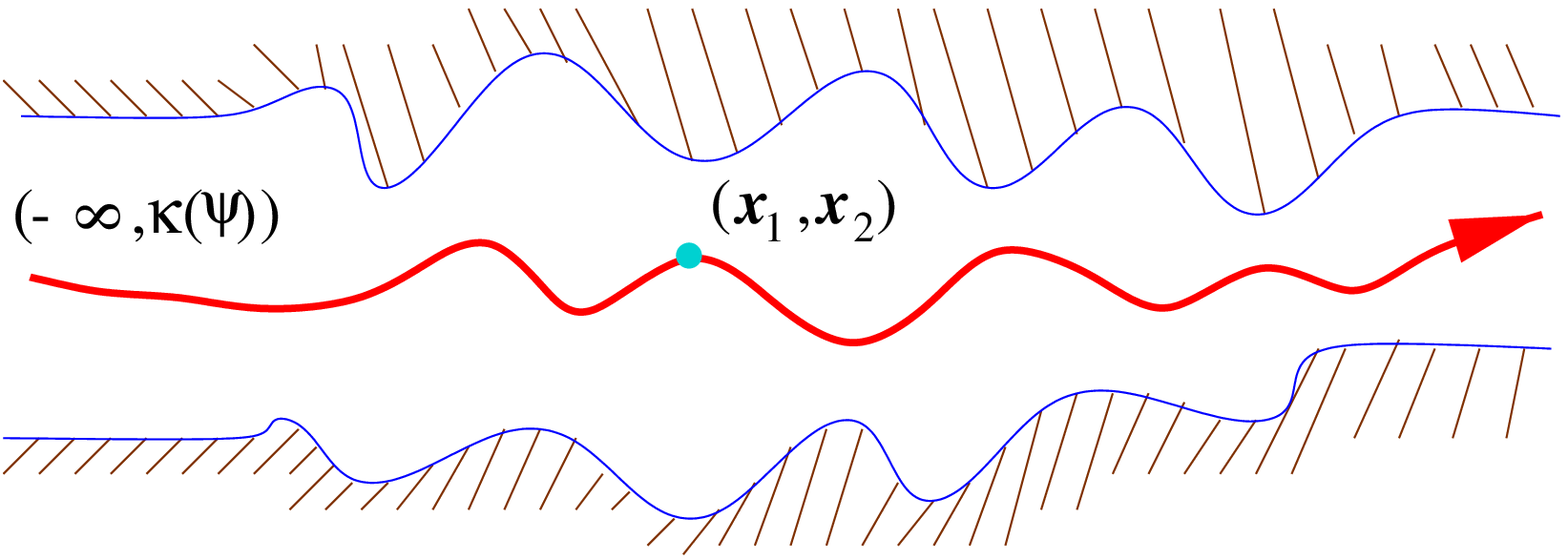}\\
     {\small The parametrization of flows by a stream function}
\end{center}

It follows from the proof of Proposition \ref{Elemmaequivalent}
that through each point $(x_1,x_2)\in\Omega$, there is one and
only one streamline which starts from the entry, provided that
(\ref{Ephorizontalvelocity}) holds in $\Omega$. By the definition
of streamlines, along each streamline, the stream function is a
constant, therefore, through any $(x_1,x_2)$ in the nozzle, there
exists a unique streamline originated from
$(-\infty,\kappa(\psi))$ with $\psi=\psi(x_1,x_2)$. Since
Bernoulli's function is also invariant along a streamline,
\begin{equation*}
(h(\rho)+\frac{|\nabla\psi|^2}{2\rho^2})(x_1,x_2)=(h(\rho)
+\frac{u^2+v^2}{2})(-\infty,\kappa(\psi))=h(\rho_0)+\frac{F^2(\psi(x_1,x_2))}{2}.
\end{equation*}
Thus in the nozzle, one has
\begin{equation}\label{EBdensitymassflux}
\mathcal{H}(\rho,|\nabla\psi|^2,\psi)=
h(\rho)+\frac{|\nabla\psi|^2}{2\rho^2}-h(\rho_0)-\frac{F^2(\psi)}{2}=0.
\end{equation}
Similarly, by virtue of (\ref{Evorticitystreamrelation}), one has
\begin{equation}\label{Eglobalstreamvorticity}
\frac{\omega}{\rho}(x_1,x_2)=-\frac{f(\psi(x_1,x_2))}{\rho_0},
\end{equation}
provided that (\ref{Ephorizontalvelocity}) holds. Furthermore,
note that (\ref{Ephorizontalvelocity}) implies
\begin{equation}\label{Estreambasicestimate}
0\leq \psi\leq m.
\end{equation}
Thus both (\ref{EBdensitymassflux}) and
(\ref{Eglobalstreamvorticity}) do make sense.

Next, we study the relationship between $F$ and $f$. In the
upstream,
\begin{equation*}
\psi=\int_0^{\kappa(\psi)}\rho_0u_0(s)ds,
\end{equation*}
which yields
\begin{equation*}
\kappa'(\psi)=\frac{1}{\rho_0u_0(\kappa(\psi))}=\frac{1}{\rho_0F(\psi)}.
\end{equation*}
So, (\ref{Eupstreamviastream}) shows
\begin{equation*}
F'(\psi)=u_0'(\kappa(\psi))\kappa'(\psi)=f(\psi)\frac{1}{\rho_0F(\psi)},
\end{equation*}
this implies
\begin{equation}
f(\psi)=\rho_0F(\psi)F'(\psi).
\end{equation}
Furthermore, if $B$ satisfies (\ref{EassumptiononBernoulli}) with
$0\leq \delta\leq \bar{\delta}_0$, and $m\in
(\delta^{\gamma},\bar{m})$, then
\begin{equation}\label{Ebstreamvariableestimate}
\left\{
\begin{array}{ll}
C^{-1}\delta^{\gamma}\leq F\leq C,\\
F'(m)\geq 0\,\, \text{and}\,\,F'(0)\leq 0,\\
|F'(\psi)|=\left|\frac{u_0'(\kappa(\psi))}{\rho_0u_0(\kappa(\psi))}\right|\leq
C\delta^{1-2\gamma},\\[1mm]
[F'(\psi)]_{C^{0,1}([0,m])} \leq  C\delta^{1-3\gamma}.
\end{array}
\right.
\end{equation}

It follows from (\ref{Edefsubsonicbranch}) and
(\ref{EBdensitymassflux}) that the subsonic flows in the nozzle
satisfy
\begin{equation}
\rho=J(|\nabla\psi|^2,h(\rho_0)+\frac{F^2(\psi)}{2}),
\end{equation}
if they have asymptotic behavior (\ref{Easymptoticupstream1}).
Furthermore, by the definitions of vorticity and stream function,
one has
\begin{equation*}
\omega=-\Div\left(\frac{\nabla\psi}{\rho}\right).
\end{equation*}
Thus, the stream function satisfies
\begin{equation}\label{Estreameq}
\Div \left(\frac{\nabla
\psi}{H(|\nabla\psi|^2,\psi)}\right)=F(\psi)F'(\psi)H(|\nabla\psi|^2,\psi),
\end{equation}
where
\begin{equation}
H(|\nabla\psi|^2,\psi)=J(|\nabla\psi|^2,h(\rho_0)+\frac{F^2(\psi)}{2}).
\end{equation}

Our major task in the rest of the paper is to show the existence
of solutions to the following boundary value problem
\begin{equation}\label{Estreampb}
\left\{
\begin{array}{ll}
\Div \left(\frac{\nabla
\psi}{H(|\nabla\psi|^2,\psi)}\right)=F(\psi)F'(\psi)H(|\nabla\psi|^2,\psi)
\,\,\text{in}\,\,\Omega,\\
\psi=\frac{x_2-f_1(x_1)}{f_2(x_1)-f_1(x_1)}m\,\,\text{on}\,\,\partial\Omega,
\end{array}
\right.
\end{equation}
and show that the flow field induced by
\begin{equation*}
\rho=H(|\nabla\psi|^2,\psi),\,\,
u=\frac{\psi_{x_2}}{\rho},\,\,v=-\frac{\psi_{x_1}}{\rho}
\end{equation*}
satisfies
(\ref{Easymptoticupstream1})-(\ref{Easymptoticdownstream2}). We
will obtain further the estimates (\ref{Estreambasicestimate}) and
(\ref{Ephorizontalvelocity}) for the solution to (\ref{Estreampb})
in order to get the existence of Euler flows.

\begin{rmk}
{\rm It is easy to see that if $B=\underline{B}$, i.e., the flow
has uniform Bernoulli's constant, then, $F$ is a constant and the
equation (\ref{Estreameq}) reduces to
\begin{equation*}
\Div\left(\frac{\nabla\psi}{H(|\nabla\psi|^2)}\right)=0.
\end{equation*}
This is nothing but the potential equation. Therefore, it is
reasonable to use (\ref{Estreampb}) to formulate the problem for
the Euler flows through the nozzles as a perturbation of the
potential flows. The condition $m>\delta^{\gamma}(\gamma<1/3)$
ensures
\begin{equation*}
|F|\geq C^{-1}\delta^{\gamma} \text{ and } |FF'|\leq
\delta^{1-2\gamma},
\end{equation*}
which guarantees that the magnitude of vorticity $|FF'|$ is
sufficiently small, thus one can regard the potential flow as a
leading ansatz for the Euler flow.}
\end{rmk}

\section{Existence of a Modified Boundary Value Problem for Stream
Function}\label{SEexistence}

There are two main difficulties to solve the problem
(\ref{Estreampb}).  The first difficulty is that the equation in
(\ref{Estreampb}) may become degenerate elliptic at sonic states. In
addition, $H$ is not well-defined for arbitrary $\psi$ and
$|\nabla\psi|$. The second difficulty  is that this is a problem in
an unbounded domain. Our basic strategy is that we extend the
definition of $F$ appropriately, truncate $|\nabla\psi|$ appeared in
$H$ in a suitable way, and use a sequence of problems on bounded
domains to approximate the orginal problem. In this section we first
get the existence of a modified problem on the unbounded domain,
which indeed solves the original problem together with the
asymptotic behavior established in the next section.

Set
\begin{equation*}
\tilde{g}(s)=\left\{
\begin{array}{ll}
F'(s),\,\,&\text{if}\,\,0\leq s\leq m,\\
F'(m)(2m-s)/m,\,\, &\text{if}\,\, m\leq s\leq 2m,\\
F'(0)(s+m)/m,\,\, &\text{if}\,\, -m\leq s\leq 0,\\
0,\,\,&\text{if}\,\, s\geq 2m,\,\,\text{or}\,\, s\leq -m. \\
\end{array}
\right.
\end{equation*}
It is obvious that $\tilde{g}\in C^{0,1}(\mathbb{R})$ and
\begin{equation*}
\|\tilde{g}(s)\|_{C^{0}(\mathbb{R}^1)}\leq \|F'(s)\|_{C^{0}([0,m])}.
\end{equation*}
Moreover, it follows from (\ref{Ebstreamvariableestimate}) that
\begin{equation}
\tilde{g}(s)\geq 0\,\,\text{if}\,\,s\geq m\,\,\text{and}\,\,
\tilde{g}(s)\leq 0\,\,\text{if}\,\,s\leq 0.
\end{equation}
Furthermore, it follows from $\|F'(s)\|_{C^{0}([0,m])}\leq
C\delta^{1-2\gamma}$, $\|F'(s)\|_{C^{0,1}([0,m])}\leq
C\delta^{1-3\gamma}$ and $m>\delta^{\gamma}$, that
\begin{equation}\label{Rextest}
\|\tilde{g}(s)\|_{C^{0,1}(\mathbb{R}^1)}\leq C\delta^{1-3\gamma}.
\end{equation}

Define
\begin{equation*}
\tilde{F}(s)=F(0)+\int_0^s\tilde{g}(t)dt.
\end{equation*}
Then $\tilde{F}'=\tilde{g}$ and $\tilde{F}\in C^{1,1}(\mathbb{R})$.
Moreover, because $m>\delta^{\gamma}$, there exists a suitably small
$\bar{\delta}_1$ such that when $\delta<\bar{\delta}_1$,
\begin{equation}
B_0<\underline{B}-\varepsilon_0\leq
h(\rho_0)+\frac{\tilde{F}^2(s)}{2}\leq \bar{B}+\varepsilon_0
\end{equation}
holds for some $\varepsilon_0>0$, where
$\bar{B}=\sup_{x_2\in[0,1]}B(x_2)$. Moreover,
(\ref{Ebstreamvariableestimate}) and (\ref{Rextest}) imply
\begin{equation}\label{Eextendfarestimatebystream}
\|\tilde{F}'\|_{C^{0}(\mathbb{R}^1)} \leq
C\delta^{1-2\gamma}\,\,\,\,\text{and}\,\,\,\,
\|\tilde{F}'\|_{C^{0,1}(\mathbb{R}^1)} \leq C\delta^{1-3\gamma}.
\end{equation}

In the rest of the paper, we will always use the following notations
\begin{equation*}
H_1(|\nabla\psi|^2,\psi)=\frac{\partial H}{\partial
|\nabla\psi|^2}(|\nabla\psi|^2,\psi),
\,\,H_2(|\nabla\psi|^2,\psi)=\frac{\partial H}{\partial
\psi}(|\nabla\psi|^2,\psi).
\end{equation*}
It follows from direct calculations that
\begin{equation*}
H_1(|\nabla\psi|^2,\psi)=-\frac{1}{2\rho(c^2-\frac{|\nabla\psi|^2}{\rho^2})}
\end{equation*}
may go to negative infinity when the flow approaches sonic from
subsonic.

Choose a smooth increasing function $\zeta_0$ such that
\begin{equation*}
\zeta_0(s)=\left\{
\begin{array}{ll}
s,\,\, &\text{if}\,\, s<-2\varepsilon_0,\\
-\varepsilon_0,\,\,&\text{if}\,\,s\geq -\varepsilon_0.
\end{array}
\right.
\end{equation*}
Then define
\begin{equation}
\tilde{\Delta}
(|\nabla\psi|^2,\psi)=\zeta_0(|\nabla\psi|^2-\Sigma^2(\tilde{\mathcal{B}}(\psi)))
+\Sigma^2(\tilde{\mathcal{B}}(\psi)),
\end{equation}
where $\Sigma$ is the function defined in
(\ref{Edefcirticalmassfluxfunction}) and
\begin{equation}
\tilde{\mathcal{B}}(\psi)=h(\rho_0)+\frac{\tilde{F}^2(\psi)}{2}.
\end{equation}
Set
\begin{equation}\label{Edeftruncateddensity}
\tilde{H}(|\nabla\psi|^2,\psi)=J(\tilde{\Delta}(|\nabla\psi|^2,\psi),
h(\rho_0)+\frac{\tilde{F}^2(\psi)}{2}),
\end{equation}
where $J$ is the function defined in (\ref{Edefsubsonicbranch}). A
direct calculation shows
\begin{equation*}
\tilde{H}_1(|\nabla\psi|^2,\psi)=-\frac{\zeta_0'\tilde{H}}
{2(\tilde{H}^2c^2-\tilde{\Delta})}.
\end{equation*}
Obviously, there exist two positive constants
$\lambda(\varepsilon_0)$ and $\Lambda(\varepsilon_0)$ such that
\begin{equation}\label{Estreamcutoffbddpb}
\lambda|\xi|^2\leq \tilde{A}_{ij}(q,z) \xi_i\xi_j \leq
\Lambda|\xi|^2
\end{equation}
holds for any $z\in\mathbb{R}^1$, $q\in \mathbb{R}^2$ and
$\xi\in\mathbb{R}^2$, where
\begin{equation}\label{Esecondordercoefficients}
 \tilde{A}_{ij}(q,z)=\tilde{H}(|q|^2,z)\delta_{ij}-
2\tilde{H}_1(|q|^2,z)q_i q_j.
\end{equation}

Instead of (\ref{Estreampb}), we first solve the following problem
\begin{equation}\label{Ecutoffstreampb}
\left\{
\begin{array}{ll}
\Div\left(\frac{\nabla\psi}{\tilde{H}(|\nabla\psi|^2,\psi)}\right)
=\tilde{F}(\psi)\tilde{F}'(\psi)\tilde{H}(|\nabla\psi|^2,\psi)\,\,\text{in}\,\,\Omega,\\
\psi=\frac{x_2-f_1(x_1)}{f_2(x_1)-f_1(x_1)}m\,\,\text{on}\,\,\partial\Omega.
\end{array}
\right.
\end{equation}

\begin{prop}\label{Elemmaexistence}
Let the boundary of $\Omega$ satisfy
(\ref{Enontrivialbdy})-(\ref{Eboundary3}). Then there exists
$0<\delta_1\leq\min\{\bar{\delta}_0,\bar{\delta}_1\}$, where
$\bar{\delta}_0$ is defined in Section \ref{SEreformulation}, such
that if $\|B'\|_{C^{0,1}([0,1])}=\delta \leq \delta_1$ and $m\in
(\delta^{\gamma},m_1)$ with $m_1=2\delta_1^{\gamma/2}\leq
\bar{m}$, where $\bar{m}$ is defined in (\ref{Edefmaxflux}) in
Section \ref{SEreformulation}, then the problem
(\ref{Ecutoffstreampb}) has a solution $\psi\in
C^{2,\alpha}(\bar{\Omega})$ satisfying
\begin{equation}\label{Estreamgradientestimate}
|\psi|\leq C(\varepsilon_0,\delta),\,\, |\nabla\psi|^2\leq
\Sigma^2(\underline{B}-\varepsilon_0)-2\varepsilon_0 \,\,\text{for
some}\,\,\varepsilon_0>0.
\end{equation}
\end{prop}
\begin{pf}
Note that the equation (\ref{Ecutoffstreampb}) is uniformly
elliptic and the domain is unbounded, one can use a sequence of
boundary value problems on bounded domains to approximate it. The
key point is to obtain the estimate
(\ref{Estreamgradientestimate}). Therefore, we first solve the
following boundary value problem
\begin{equation}\label{Etruncatedbddpb}
\left\{
\begin{array}{ll}
\Div\left(\frac{\nabla\psi}{\tilde{H}(|\nabla\psi|^2,\psi)}\right)
=\tilde{F}(\psi)\tilde{F}'(\psi)\tilde{H}(|\nabla\psi|^2,\psi)\,\, \text{in}\,\,\Omega_L,\\
\psi=\frac{x_2-f_1(x_1)}{f_2(x_1)-f_1(x_1)}m\,\,\text{on}\,\,\partial\Omega_L,
\end{array}
\right.
\end{equation}
where $\Omega_L$ satisfies $\{(x_1,x_2)|(x_1,x_2)\in\Omega,
-L<x_1<L\}\subset\Omega_L\subset\{(x_1,x_2)|(x_1,x_2)\in\Omega,
-4L<x_1<4L\}$ for $\forall L\in\mathbb{N}$. Furthermore, one may
choose $\Omega_L$ so that $\Omega_L\in
C^{2,\alpha_1}$($0<\alpha_1\leq \alpha$) satisfies the uniform
exterior sphere condition with uniform radius $r_0$, $0<r_0<r$, for
all $L>L_0$ with some $L_0$ sufficiently large. For the explicit
construction of such $\Omega_L$, please refer to Appendix in
\cite{XX1}.

The equation in (\ref{Ecutoffstreampb}) can be rewritten as
\begin{equation}\label{Etruncatednondiveq1}
\tilde{A}_{ij}(D\psi,\psi)\partial_{ij}\psi
-\tilde{H}_2(|\nabla\psi|^2,\psi)|\nabla\psi|^2
=\tilde{F}(\psi)\tilde{F}'(\psi)\tilde{H}^3(|\nabla\psi|^2,\psi),
\end{equation}
where
\begin{equation*}
\tilde{H}_2(|\nabla\psi|^2,\psi)=\frac{\tilde{F}(\psi)\tilde{F}'(\psi)
\tilde{H}(|\nabla\psi|^2,\psi)(\tilde{H}^2+\Sigma\Sigma'(\zeta_0'-1))}
{\tilde{H}^2(|\nabla\psi|^2,\psi)c^2
-\tilde{\Delta}(|\nabla\psi|^2,\psi)}.
\end{equation*}
Therefore, (\ref{Etruncatednondiveq1}) becomes
\begin{equation}\label{Etruncatednondiveqcpt}
\tilde{A}_{ij}(D\psi,\psi)\partial_{ij}\psi
=\mathcal{F}(\psi,\nabla\psi),
\end{equation}
where
\begin{equation*}
\mathcal{F}(\psi,\nabla\psi)=\tilde{F}(\psi)\tilde{F}'(\psi)\tilde{H}(|\nabla\psi|^2,\psi)
\left(\frac{(\tilde{H}^2+\Sigma\Sigma'(\zeta_0'-1))|\nabla\psi|^2}
{\tilde{H}^2(|\nabla\psi|^2,\psi)c^2
-\tilde{\Delta}(|\nabla\psi|^2,\psi)} +\tilde{H}^2\right).
\end{equation*}
Note that $\mathcal{F}$ has quadratic growth in $|\nabla\psi|$, so
it is not easy to get a prior estimate and the existence for
(\ref{Etruncatednondiveqcpt}) directly. The strategy here is that,
instead of (\ref{Etruncatedbddpb}), we first solve the problem
\begin{equation}\label{Etruncatedbddcutoffpb}
\left\{
\begin{array}{ll}
\tilde{A}_{ij}(D\psi,\psi)\partial_{ij}\psi
=\tilde{\mathcal{F}}(\psi,\nabla\psi)\,\,\text{in}\,\,\Omega_L,\\
\psi=\frac{x_2-f_1(x_1)}{f_2(x_1)-f_1(x_1)}m\,\,\text{on}\,\,\partial\Omega_L,
\end{array}
\right.
\end{equation}
where
\begin{equation*}
\tilde{\mathcal{F}}(\psi,\nabla\psi)=\tilde{F}(\psi)\tilde{F}'(\psi)\tilde{H}(|\nabla\psi|^2,\psi)
\left(\frac{(\tilde{H}^2+\Sigma\Sigma'(\zeta_0'-1))\tilde{\Delta}(\nabla\psi,\psi)}
{\tilde{H}^2(|\nabla\psi|^2,\psi)c^2
-\tilde{\Delta}(|\nabla\psi|^2,\psi)} +\tilde{H}^2\right).
\end{equation*}

Thanks to (\ref{Eextendfarestimatebystream}), one has
\begin{equation}\label{ERHSestimate}
|\tilde{\mathcal{F}}|\leq C\delta^{1-2\gamma}.
\end{equation}
It follows from Theorem 12.5 and Remark in P308 in \cite{GT} that
there exists a solution $\psi_L$ to (\ref{Etruncatedbddcutoffpb}).
Furthermore, writing $\psi_L^-=\min\{\psi_L,0\}$ and
$\psi_L^{+}=\max\{\psi_L,0\}$, by the proof of Theorem 3.7 in
\cite{GT},
\begin{equation}\label{Eaprioribdcutoffpb}
\inf_{\partial\Omega_L}\psi_L^{-}-C\sup_{\Omega_L}
\left|\frac{\tilde{\mathcal{F}}}{\lambda}\right|\leq \psi_L\leq
\sup_{\partial\Omega_L}\psi_L^{+}+C\sup_{\Omega_L}\left|\frac{\tilde{\mathcal{F}}}
{\lambda}\right|,
\end{equation}
where $C=e^d-1$ with $d=\sup\{f_2(x_1)-f_1(x_1)\}$. Thus,
\begin{equation*}
-C\delta^{1-2\gamma}\leq \psi_k\leq m+C\delta^{1-2\gamma} \,\,
\text{for}\,\, k\,\, \text{sufficiently large}.
\end{equation*}

Moreover, one can get some nice estimates for $\psi_k$. This follows
from the techniques developed in Chapter 12 in \cite{GT}. Using the
specific form of estimate (12.14) in P299 in \cite{GT} and Remark
(4) on global estimate for quasiconformal mappings in P300 in
\cite{GT}, then one can improve the estimate in Line 7 in P304 in
\cite{GT} to the following more precise form
\begin{equation}\label{EimproveGTestimate}
[u]_{1,\alpha}\leq
C(\gamma,\Omega)\left(1+|Du|_0+\left|\frac{f}{\lambda}\right|_0\right),
\end{equation}
actually, $C(\gamma,\Omega)$ depends only on the diam$\Omega$ and
$C^2$ norm of $\partial\Omega$. Here we use notations and symbols in
(\ref{EimproveGTestimate}) as those in Chapter 12 in \cite{GT}.

Note that although the estimate (\ref{EimproveGTestimate}) is
derived with zero boundary conditions, it holds in the case that
the boundary value is constant in each connect component of
boundary. Indeed, first, it holds for the case that the boundary
value is a constant. Then one can generalize the estimate to the
case that boundary value is constant in each connected component
of the boundary, since all estimates are obtained through
localization.

Applying the estimate (\ref{EimproveGTestimate}) to the problem
(\ref{Etruncatedbddcutoffpb}) shows that, there exists
$\mu=\mu(\Lambda/\lambda)>0$, such that for any $x^0\in
\bar{\Omega}_L$, and for $\psi_k$ with $k\geq 4L$, one has
\begin{equation}\label{EbasicHestimate}
[\psi_k]_{1,\mu;B_1(x^0)\cap\Omega_L}\leq
C(\Lambda/\lambda,|f_i|_2)\left(1+|D\psi_k|_{0;B_1(x^0)\cap\Omega_L}+\left|\frac{\tilde{\mathcal{F}}}
{\lambda}\right|_0\right).
\end{equation}
This, together with interpolation inequality and
(\ref{Eaprioribdcutoffpb}), yields
\begin{align*}
\|\psi_k\|_{1;B_1(x^0)\cap\Omega_L} \leq & \eta
[\psi_k]_{1,\mu;B_1(x^0)\cap\Omega_L}+C_{\eta}|\psi_k|_0\nonumber \\
\leq & \eta C(\Lambda/\lambda,|f_i|_2)
\left(1+|D\psi_k|_{0;B_1(x^0)\cap\Omega_L}+
\left|\frac{\tilde{\mathcal{F}}}{\lambda}\right|_0\right)
+C_{\eta}\left(m+C\left|\frac{\tilde{\mathcal{F}}}{\lambda}\right|_0\right),
\end{align*}
where $C$, appeared in last term, is the same as that in
(\ref{Eaprioribdcutoffpb}).  Taking $\eta_0$ sufficiently small so
that $\eta C(\Lambda/\lambda,|f_i|_2)\leq 1/2$ if $\eta\leq\eta_0$,
then one has
\begin{equation}\label{Epgestimate}
\|\psi_k\|_{1;B_1(x^0)\cap\Omega_L} \leq \eta
C(\Lambda/\lambda,|f_i|_2)
\left(1+\left|\frac{\tilde{\mathcal{F}}}{\lambda}\right|_0\right)
+C_{\eta}\left(m+C\left|\frac{\tilde{\mathcal{F}}}{\lambda}\right|_0\right).
\end{equation}
Thus, the H\"older estimate (\ref{EbasicHestimate}) becomes
\begin{align}
\|\psi_k\|_{1,\mu,B_1(x^0)\cap \Omega_L}\leq &
\|\psi_k\|_{1;B_1(x^0)\cap\Omega_L}+[\psi_k]_{1,\mu;B_1(x^0)\cap\Omega_L}\nonumber\\
\leq &
(1+C(\Lambda/\lambda,|f_i|_2))\|\psi_k\|_{1;B_1(x^0)\cap\Omega_L}
+C(\Lambda/\lambda,|f_i|_2)
\left(1+\left|\frac{\tilde{\mathcal{F}}}{\lambda}\right|_0\right)\nonumber\\
\leq & C(\Lambda/\lambda,|f_i|_2)\left(\eta_0
C(\Lambda/\lambda,|f_i|_2)
\left(1+\left|\frac{\tilde{\mathcal{F}}}{\lambda}\right|_0\right)
+C_{\eta_0}\left(m+C\left|\frac{\tilde{\mathcal{F}}}{\lambda}\right|_0\right)\right)\nonumber\\
&+C(\Lambda/\lambda,|f_i|_2)
\left(1+\left|\frac{\tilde{\mathcal{F}}}{\lambda}\right|_0\right)\nonumber\\
\leq &C(\Lambda/\lambda,|f_i|_2))
\left(1+m+\left|\frac{\tilde{\mathcal{F}}}{\lambda}\right|_0\right).
\label{Elocalholderforcutpb}
\end{align}
Note that, for any $x$, $y\in \bar{\Omega}_L$,
\begin{equation*}
\frac{|\nabla\psi_k(x)-\nabla\psi_k(y)|} {|x-y|^{\mu}}\leq\left\{
\begin{aligned}
&\|\psi_k\|_{1,\mu;B_1(x)\cap\Omega_L},\,\,\text{if}\,\,y\in
B_1(x)\cap\bar{\Omega}_L,\\
&2\|\psi_k\|_{1;\Omega_L},\,\,\,\,\,\qquad\text{if}\,\,y\notin
B_1(x)\cap\bar{\Omega}_L.
\end{aligned}
\right.
\end{equation*}
This, together with (\ref{Epgestimate}) and
(\ref{Elocalholderforcutpb}), yields the following H\"older estimate
\begin{equation}\label{Eggholderestimate}
[\psi_k]_{1,\mu;\Omega_L}=\sup_{x,y\in\Omega_L}
\frac{|\nabla\psi_k(x)-\nabla\psi_k(y)|} {|x-y|^{\mu}} \leq
C(\Lambda/\lambda,|f_i|_2)
\left(1+m+\left|\frac{\tilde{\mathcal{F}}}{\lambda}\right|_0\right).
\end{equation}
Furthermore, it follows from (\ref{Elocalholderforcutpb}), the
Schauder estimate (Theorem 6.2 and Lemma 6.5 in \cite{GT}), and the
bootstrap argument that
\begin{equation*}
\|\psi_k\|_{2,\alpha; B_{1/2}(x^0)\cap\Omega_L}\leq C\left(
\Lambda/\lambda,|f_i|_{C^{2,\alpha}},m,
\left|\frac{\mathcal{\tilde{F}}}{\lambda}\right|_0\right).
\end{equation*}
Similar to the argument for (\ref{Eggholderestimate}), one has
\begin{equation}\label{ESchauderestimate}
\|\psi_k\|_{2,\alpha;\Omega_L}\leq C\left(
\Lambda/\lambda,|f_i|_{C^{2,\alpha}},m,
\left|\frac{\mathcal{\tilde{F}}}{\lambda}\right|_0\right).
\end{equation}
Hence, using Arzela-Ascoli lemma and a diagonal procedure, we see
that there exists a sequence $\psi_{k_l}$ such that
\begin{equation*}
\psi_{k_l}\rightarrow \psi\,\,\text{in} \,\,
C^{2,\beta}(K)\,\,\text{for any compact set}\,\, K\subset
\bar{\Omega}\,\, \text{and}\,\,\beta<\alpha.
\end{equation*}
Furthermore, $\psi$ satisfies the problem
\begin{equation*}
\left\{
\begin{array}{ll}
\tilde{A}_{ij}(D\psi,\psi)\partial_{ij}\psi
=\tilde{\mathcal{F}}(\nabla\psi,\psi)\,\,\text{in}\,\,\Omega,\\
\psi=\frac{x_2-f_1(x_1)}{f_2(x_1)-f_1(x_1)}m\,\,\text{on}\,\,\partial\Omega,
\end{array}
\right.
\end{equation*}
and the estimate
\begin{equation*}
\|\psi\|_{1,\Omega}\leq \eta
C(\gamma,|f_i|_2)\left(1+\left|\frac{\tilde{\mathcal{F}}}{\lambda}\right|_0\right)
+C_{\eta}\left(m+C\left|\frac{\tilde{\mathcal{F}}}{\lambda}\right|_0\right),
\end{equation*}
where $\eta \in (0,\eta_0)$. Thanks to estimate
(\ref{ERHSestimate}), one has
\begin{equation}\label{EC1estimateformpb}
\|\psi\|_{1,\Omega}\leq \eta
C(\gamma,|f_i|_2)(1+C\delta^{1-2\gamma})
+C_{\eta}(m+C\delta^{1-2\gamma}),
\end{equation}
where $C$ depends only on $\bar{\delta}_0$, $\bar{m}$, $\Lambda$ and
$\lambda$.

Obviously, there exist $\eta_1\in (0,\eta_0)$ and $\delta_1\in
(0,\bar{\delta}_0]$ such that
\begin{eqnarray*}
&\eta_1 C(\gamma,|f_i|_2)(1+C\bar{\delta}_0^{1-2\gamma})\leq
\sqrt{(\Sigma^2(\underline{B}-\varepsilon_0)-2\varepsilon_0)/2},\\
&C_{\eta_1}(2\delta_1^{\gamma/2}+C\delta_1^{1-2\gamma})\leq
\sqrt{(\Sigma^2(\underline{B}-\varepsilon_0)-2\varepsilon_0)/2}.
\end{eqnarray*}
Therefore, for any $\delta\in (0,\delta_1)$ and $m\in
(\delta^{\gamma},2\delta_1^{\gamma/2})$, the solution $\psi$
satisfies
\begin{equation}\label{Egradientformpb}
|\nabla\psi|^2\leq
\Sigma^2(\underline{B}-\varepsilon_0)-2\varepsilon_0.
\end{equation}
Now (\ref{Estreamgradientestimate}) follows from
(\ref{EC1estimateformpb}) and (\ref{Egradientformpb}).

Furthermore,  (\ref{Eggholderestimate}) and
(\ref{ESchauderestimate}) yield the following higher order estimates
\begin{align}
&\|\psi\|_{1,\mu;\bar{\Omega}}\leq
C(\Lambda/\lambda,|f_i|_2)\left(1+m+\left|\frac{\mathcal{F}}{\lambda}\right|_0
\right),\label{Eggholderestimateformpb}\\
&\|\psi\|_{2,\alpha;\bar{\Omega}}\leq C\left(\Lambda/\lambda,
|f_i|_{C^{2,\alpha}},m,
\left|\frac{\mathcal{F}}{\lambda}\right|_0\right)\label{ESchauderformpb}.
\end{align}
This finishes the proof of the Proposition.
\end{pf}

\section{Far Fields Behavior, Existence
and Uniqueness of Boundary Value Problem for the Stream
Function}\label{SEasymptotic}

In this section, we will study far fields behavior of the solution
to (\ref{Ecutoffstreampb}). It will be shown that the flows induced
by the solutions to (\ref{Ecutoffstreampb}) satisfy asymptotic
behavior (\ref{Easymptoticupstream1})-(\ref{Easymptoticupstream2}).
This also yields  that solutions to (\ref{Ecutoffstreampb}) satisfy
(\ref{Estreambasicestimate}). Combining (\ref{Estreambasicestimate})
and (\ref{Egradientformpb}), we can remove both extension and
truncation appeared in (\ref{Ecutoffstreampb}). Therefore, these
solutions solve problem (\ref{Estreampb}). Furthermore, the
asymptotic behavior is crucial for our formulation for the problem
since the stream function formulation is consistent with the
original formulation of the problem for the Euler system in the
infinitely long nozzle, as long as the flow induced by a solution to
(\ref{Estreampb}) satisfies
(\ref{Easymptoticupstream1})-(\ref{Easymptoticupstream2}) and
(\ref{Ephorizontalvelocity}). Finally, the uniqueness of the
solutions will be a consequence of the asymptotic behavior. To study
the solution in its far fields, we will use a blow up argument and
an energy estimate.

For $x_1\leq n$, define
$\psi^{(n)}(x_1,x_2)=\psi(x_1-n,x_2)\chi_{\{f_1(x_1-n)<x_2<f_2(x_1-n)\}}$.
For any compact set $K\Subset (-\infty,\infty)\times(0,1)$, it
follows from (\ref{ESchauderformpb}) that
\begin{equation*}
\|\psi^{(n)}\|_{C^{2,\alpha}(K)}\leq C\,\,\text{for}\,\,
n\,\,\text{sufficiently large.}
\end{equation*}
Therefore, by Arzela-Ascoli lemma and a diagonal procedure, there
exists a subsequence, $\psi^{(n_k)}$, such that
\begin{equation}\label{Eblowup2ndHolderconvergence}
\psi^{(n_k)}\rightarrow \bar{\psi} \,\, \text{in}\,\,C^{2,\beta}(K)
\end{equation}
for any $K\Subset(-\infty,\infty)\times(0,1)$, for any
$\beta\in(0,\alpha)$. Furthermore, it follows from
(\ref{Enontrivialbdy})-(\ref{Eboundary3}) and
(\ref{ESchauderformpb}) that $\bar{\psi}=0$ on $x_2=0$ and
$\bar{\psi}=m$ on $x_2=1$. So, $\bar{\psi}$ satisfies
\begin{equation}\label{Eblowuppb}
\left\{
\begin{array}{ll}
\Div\left(\frac{\nabla\bar{\psi}}{\tilde{H}(|\nabla\bar{\psi}|^2,\bar{\psi})}\right)
=\tilde{F}'(\bar{\psi})\tilde{F}(\bar{\psi})\tilde{H}(|\nabla\bar{\psi}|^2,\bar{\psi})\,\,\text{in}\,\, D,\\
\bar{\psi}=0\,\,\text{on}\,\, x_2=0,\,\,
\bar{\psi}=m\,\,\text{on}\,\, x_2=1,
\end{array}
\right.
\end{equation}
where $D=(-\infty,\infty)\times (0,1)$. Moreover, by
(\ref{Estreamgradientestimate}), one has
\begin{equation}\label{Eblowupc1estimate}
|\bar{\psi}|\leq C(\varepsilon_0,\delta)\,\,\text{and}\,\,
|\nabla\bar{\psi}|^2\leq
\Sigma^2(\underline{B}-\varepsilon_0)-2\varepsilon_0.
\end{equation}
Thus, by the similar argument in Section \ref{SEexistence}, on any
compact set $E\subset (-\infty,\infty)\times[0,1]$,
\begin{equation*}
\|\bar{\psi}\|_{C^{1,\mu}(E)}\leq C(\varepsilon,\delta).
\end{equation*}
Moreover, it follows from the Schauder estimate for second order
uniformly elliptic equations that
\begin{equation}\label{Eblowup2ndHolderestimate}
\|\bar{\psi}\|_{C^{2,\alpha}(E)}\leq C(\varepsilon,\delta).
\end{equation}
Therefore, $\bar{\psi}\in C^{2,\alpha}(\bar{D})$. In fact, we have
the following stronger results
\begin{lem}\label{Elemmaasymptotic}
There exists $\delta_2\in (0,\bar{\delta}_0]$ such that if
\begin{enumerate}
\item[(i).] $\|B'\|_{C^{0,1}([0,1])}=\delta\leq \delta_2$,
\item[(ii).] $m\in (\delta^{\gamma},\bar{m})$,  where
$\bar{m}$ is defined in (\ref{Edefmaxflux}) in Section
\ref{SEreformulation},
\item[(iii).] there exists $\epsilon\leq \varepsilon_0$ such that  $\bar{\psi}$ satisfies
\begin{equation}\label{Econditionforasymptotics}
|\bar{\psi}|\leq C(\epsilon,\delta)\,\,\text{and}\,\,
|\nabla\bar{\psi}|^2-\Sigma^2(\tilde{\mathcal{B}}(\bar{\psi})) \leq
-\epsilon,
\end{equation}
and solves the problem (\ref{Eblowuppb}), where
$\tilde{\mathcal{B}}$ is defined in Section \ref{SEexistence},
\end{enumerate}
then $\bar{\psi}$ is independent of $x_1$, moreover,
\begin{equation}\label{Eblowupexplicitform}
\bar{\psi}(x_1,x_2)=\bar{\psi}(x_2)=\int_0^{x_2}\rho_0u_0(s)ds,
\end{equation}
 where $\rho_0$ and $u_0$ are uniquely determined by
$B$ and $m$ as in Section \ref{SEreformulation}.
\end{lem}
\begin{pf}
The proof is divided into two steps. First, it will be shown that
$\bar{\psi}$ is independent of $x_1$. Then we will prove that
$\bar{\psi}$ is of explicit form (\ref{Eblowupexplicitform}).

Step 1. Set $w=\bar{\psi}_{x_1}$. Differentiating the equation in
(\ref{Eblowuppb}) with respect to $x_1$ yields
\begin{equation}\label{Eblowupderivativeeq}
\begin{array}{ll}
\,\,\,\,\,\,\,\partial_{i}\left(\frac{\tilde{A}_{ij}(D\bar{\psi},\bar{\psi})}
{\tilde{H}^2(|\nabla\bar{\psi}|^2,\bar{\psi})}\partial_j
w\right)-\partial_i\left(\frac{\tilde{H}_2(|\nabla\bar{\psi}|^2,\bar{\psi})
\partial_i\bar{\psi}}{\tilde{H}^2(|\nabla\bar{\psi}|^2,\bar{\psi})}w\right)\\
=\tilde{\Theta}(|\nabla\bar{\psi}|^2,\bar{\psi}) w
+\tilde{\vartheta}(|\nabla\bar{\psi}|^2,\bar{\psi})
\partial_i\bar{\psi}\partial_i w,
\end{array}
\end{equation}
where $\tilde{A}_{ij}$, $\tilde{\Theta}$ and $\tilde{\vartheta}$ are
defined as
\begin{eqnarray}
\tilde{A}_{ij}(q,z)&=&\tilde{H}(|q|^2,z)\delta_{ij}-
2\tilde{H}_1(|q|^2,z)q_i q_j,\label{Eprincipalcoefficient}\\
\tilde{\Theta}(s,z)&=&(\tilde{F}^{''}(z)\tilde{F}(z)+(\tilde{F}'(z))^2)\tilde{H}(s,z)
+\tilde{F}'(z)\tilde{F}(z)\tilde{H}_2(s,z),\label{Ederivativeeqright1}\\
\tilde{\vartheta}(s,z)&=&2\tilde{F}(z)\tilde{F}'(z)\tilde{H}_1(s,z),\label{Ederivativeeqright2}
\end{eqnarray}
for any $q\in \mathbb{R}^2$, $s\geq 0$, and $z\in\mathbb{R}$, here
$\tilde{F}''\in L^{\infty}(\mathbb{R}^1)$ since $\|\tilde{F}'\|_{
C^{0,1}(\mathbb{R}^1)}\leq C\delta^{1-3\gamma}$. It follows from
(\ref{Econditionforasymptotics}) that there exists a constant
$\Lambda$ depending only on $\epsilon$ such that
\begin{equation*}
|\tilde{A}_{ij}(D\bar{\psi},\bar{\psi})|\leq\Lambda(\epsilon).
\end{equation*}
Although it is unknown whether $\bar{\psi}\in C^3(D)$, the equation
(\ref{Eblowupderivativeeq}) holds in weak sense. Moreover, $w$
satisfies the boundary conditions
\begin{equation*}
w=0\,\,\text{on}\,\,x_2=0,1.
\end{equation*}

Let $\eta$ be a $C_{0}^{\infty}$ function satisfying
\begin{equation}\label{Ecutofftestfunction}
\eta=1\,\, \text{for}\,\,|s|<l,\,\,\eta=0\,\,\text{for}\,\,|s|>l+1,
\,\,\text{and}\,\,|\eta'(s)|\leq 2.
\end{equation}
Now multiplying $\eta^2(x_1)w$ on both sides of
(\ref{Eblowupderivativeeq}) and integrating it over $D$ yield
\begin{eqnarray*}
&&\iint_D
\frac{\tilde{A}_{ij}(D\bar{\psi},\bar{\psi})}{\tilde{H}^2(|\nabla\bar{\psi}|^2,\bar{\psi})}\partial_j
w \partial_i(\eta^2 w)-
\frac{\tilde{H}_2(|\nabla\bar{\psi}|^2,\bar{\psi})\partial_i\bar{\psi}}{\tilde{H}^2(|\nabla\bar{\psi}|^2,\bar{\psi})}w\partial_i(\eta^2
w)dx_1dx_2\\
&=&-\iint_D \tilde{\Theta}(|\nabla\bar{\psi}|^2,\bar{\psi}) \eta^2
w^2 dx_1dx_2+
\tilde{\vartheta}(|\nabla\bar{\psi}|^2,\bar{\psi})\partial_i\bar{\psi}\partial_i
w \eta^2 wdx_1dx_2.
\end{eqnarray*}
Substituting the explicit forms of $A_{ij}$,
$\tilde{H}_1(|\nabla\bar{\psi}|^2,\bar{\psi})$ and
$\tilde{H}_2(|\nabla\bar{\psi}|^2,\bar{\psi})$ into the above
equality and noting that $\bar{\psi}$ satisfies
(\ref{Econditionforasymptotics}), one may get
\begin{eqnarray*}
&&\iint_D \frac{\eta^2|\nabla w|^2}{\tilde{H}(|\nabla\bar{\psi}|^2,\bar{\psi})}dx_1 dx_2\\
&=&\iint_D\frac{2\tilde{H}_1(|\nabla\bar{\psi}|^2,\bar{\psi})}{\tilde{H}^2
(|\nabla\bar{\psi}|^2,\bar{\psi})}|\nabla\psi\cdot\nabla
w|^2\eta^2 dx_1dx_2\\
&&-\iint_D
2\frac{\tilde{A}_{ij}(D\bar{\psi},\bar{\psi})}{\tilde{H}^2
(|\nabla\bar{\psi}|^2,\bar{\psi})}\partial_j w\partial_i \eta \eta wdx_1dx_2\\
&&+ \iint_D
\frac{\tilde{H}_2(|\nabla\bar{\psi}|^2,\bar{\psi})\partial_i\bar{\psi}}
{\tilde{H}^2(|\nabla\bar{\psi}|^2,\bar{\psi})}(\eta^2 w\partial_i
w+2\eta
\partial_i\eta w^2)dx_1dx_2 \\
&&-\iint_D
(\tilde{F}^{''}(\bar{\psi})\tilde{F}(\bar{\psi})+(\tilde{F}'(\bar{\psi}))^2)
\tilde{H}(|\nabla\bar{\psi}|^2,\bar{\psi})\eta^2 w^2dx_1dx_2\\
&& -\iint_D
\tilde{F}'(\bar{\psi})\tilde{F}(\bar{\psi})\tilde{H}_2(|\nabla\bar{\psi}|^2,\bar{\psi})
\eta^2 w^2dx_1dx_2\\
&&-\iint_D
2\tilde{F}'(\bar{\psi})\tilde{F}(\bar{\psi})\tilde{H}_1(|\nabla\bar{\psi}|^2,\bar{\psi})
\eta^2 w\nabla\bar{\psi}\cdot\nabla wdx_1dx_2\\
&=& - \iint_D \frac{|\nabla\bar{\psi}\cdot\nabla w|^2\eta^2}
{\tilde{H}(|\nabla\bar{\psi}|^2,\bar{\psi})(\tilde{H}^2(|\nabla\bar{\psi}|^2,\bar{\psi})c^2
-|\nabla\bar{\psi}|^2)}dx_1dx_2\\
&&-\iint_D
2\frac{\tilde{A}_{ij}(D\bar{\psi},\bar{\psi})}{\tilde{H}^2
(|\nabla\bar{\psi}|^2,\bar{\psi})}\partial_j w\partial_i \eta \eta
wdx_1dx_2\\
&& + \iint_D
\frac{2\tilde{H}_2(|\nabla\bar{\psi}|^2,\bar{\psi})\nabla\bar{\psi}\cdot\nabla\eta}
{\tilde{H}^2(|\nabla\bar{\psi}|^2,\bar{\psi})}\eta w^2 dx_1dx_2\\
&&+\iint_D
\frac{2\tilde{F}(\bar{\psi})\tilde{F}'(\bar{\psi})\tilde{H} \nabla
\bar{\psi}\cdot\nabla w}
{\tilde{H}^2(|\nabla\bar{\psi}|^2,\bar{\psi}) c^2
-|\nabla\bar{\psi}|^2}
\eta^2 wdx_1dx_2\\
&&-\iint_D \left(\tilde{F}''(\bar{\psi})\tilde{F}(\bar{\psi})+
(\tilde{F}'(\bar{\psi}))^2\right)
\tilde{H}(|\nabla\bar{\psi}|^2,\bar{\psi})
\eta^2 w^2 dx_1dx_2\\
&&-\iint_D \frac{(\tilde{F}(\psi)\tilde{F}'(\psi))^2
\tilde{H}^3(|\nabla\bar{\psi}|^2,\bar{\psi}) }
{\tilde{H}^2(|\nabla\bar{\psi}|^2,\bar{\psi})c^2
-|\nabla\bar{\psi}|^2} \eta^2 w^2dx_1dx_2,
\end{eqnarray*}
which can be written as
\begin{equation}\label{Ecptformkeyest}
\iint_D \frac{\eta^2|\nabla
w|^2}{\tilde{H}(|\nabla\bar{\psi}|^2,\bar{\psi})}dx_1
dx_2=\sum_{i=1}^6 I_i.
\end{equation}
First, it is easy to see that $I_1+I_4+I_6\leq 0$. Second, due to
(\ref{Eextendfarestimatebystream}), one has
\begin{equation}\label{EestI5}
|I_5|\leq C\delta^{1-3\gamma}\int_{-l-1}^{l+1}\int_0^1w^2
dx_1dx_2.
\end{equation}
Finally, since $\tilde{H}\leq \bar{\varrho}(\bar{B})$, thus if
$\delta_2$ is sufficiently small, one gets from above that
\begin{eqnarray*}
&&\int_{-l}^l dx_1\int_0^1|\nabla w|^2 dx_2\\
&\leq&
C(\bar{B},\epsilon)\left(\int_{-l-1}^{-l}dx_1+\int_l^{l+1}dx_1\right)\int_0^1|\nabla
w|^2 +|\nabla w w|+w^2dx_2\\
& &+C(\bar{B})\delta^{1-3\gamma}\int_{-l}^l\int_0^1w^2dx_1dx_2\\
&\leq &
C(\bar{B},\epsilon)\left(\int_{-l-1}^{-l}dx_1+\int_l^{l+1}dx_1\right)\int_0^1|\nabla
w|^2
+w^2dx_2+C(\bar{B})\delta^{1-3\gamma}\int_{-l}^l\int_0^1w^2dx_1dx_2.
\end{eqnarray*}
Notice that $w=0$ on $x_2=0$. It follows from Poincare inequality
that
\begin{equation}\label{Epoincareinequality}
\int_0^1w^2dx_2\leq \int_0^1|\nabla  w|^2dx_2.
\end{equation}
Therefore, there exists a constant $C$ independent of $l$ such that
\begin{equation}\label{Eblowupgradientbyfiniteinterval}
\int_{-l}^l \int_0^1|\nabla w|^2 dx_1dx_2 \leq
C\left(\int_{-l-1}^{-l}dx_1+\int_l^{l+1}dx_1\right)\int_0^1|\nabla
w|^2dx_2
\end{equation}
 for large $l$. It follows from
(\ref{Eblowup2ndHolderestimate}) that
\begin{equation*}
\left(\int_{-l-1}^{-l}dx_1+\int_l^{l+1}dx_1\right)\int_0^1|\nabla
w|^2dx_2\leq C
\end{equation*}
for some uniform constant $C$  independent of $l$. Therefore,
\begin{equation*}
\int_{-l}^l dx_1\int_0^1|\nabla w|^2 dx_2 \leq C
\end{equation*}
for some constant $C$. Taking  $l\rightarrow \infty$ yields
\begin{equation*}
\int_{-\infty}^{\infty} dx_1\int_0^1|\nabla w|^2 dx_2 \leq C.
\end{equation*}
Hence
\begin{equation}\label{Evanishingfarintegral}
\left(\int_{-l-1}^{-l}dx_1+\int_l^{l+1}dx_1\right)\int_0^1|\nabla
w|^2dx_2\rightarrow 0\,\,\text{as}\,\,l\rightarrow \infty.
\end{equation}
Taking the limit $l\rightarrow \infty$ in
(\ref{Eblowupgradientbyfiniteinterval}), one has
\begin{equation*}
\int_{-\infty}^{\infty} \int_0^1|\nabla w|^2 dx_1 dx_2 =0.
\end{equation*}
So
\begin{equation*}
w=0.
\end{equation*}

Therefore, $\bar{\psi}=\bar{\psi}(x_2)$. Thus $\bar{\psi}$ solves
the following boundary value problem
\begin{equation}\label{EblowupODEBVP}
\left\{
\begin{array}{ll}
\frac{d}{dx_2}\left(\frac{\bar{\psi}'}{\tilde{H}((\bar{\psi}')^2,\bar{\psi})}\right)=\tilde{F}'(\psi)
\tilde{F}(\psi)\tilde{H}((\bar{\psi}')^2,\bar{\psi}),\\
\bar{\psi}(0)=0,\,\,\bar{\psi}(1)=m.
\end{array}
\right.
\end{equation}

Step 2. Uniqueness of the solution to the boundary value problem
(\ref{EblowupODEBVP}).

Suppose that there are two solutions $\bar{\psi}_1$ and
$\bar{\psi}_2$ to (\ref{EblowupODEBVP}). Let
$\bar{\phi}=\bar{\psi}_1-\bar{\psi}_2$. Then $\bar{\phi}$ satisfies
\begin{equation}\label{EblowupODEdifference}
\left\{
\begin{array}{ll}
(\bar{a}\bar{\phi}'+\bar{b}\bar{\phi})'=
\bar{c}\bar{\phi}'+\bar{d}\bar{\phi},\\
\bar{\phi}(0)=\bar{\phi}(1)=0,
\end{array}
\right.
\end{equation}
where
\begin{eqnarray*}
\bar{a}&=&\int_0^1 \frac{\tilde{H}(|\tilde{\psi}|^2,
\tilde{\psi})-2\tilde{H}_1 (|\tilde{\psi}'|^2,
\tilde{\psi})|\tilde{\psi}'|^2
}{\tilde{H}^2(|\tilde{\psi}'|,\tilde{\psi})} ds,\,\,
\bar{b}=\int_0^1 \frac{-\tilde{H}_2 (|\tilde{\psi}'|^2,
\tilde{\psi})\tilde{\psi}'}
{\tilde{H}^2(|\tilde{\psi}'|,\tilde{\psi})} ds,\\
\bar{c}&=&\int_0^1
\tilde{\vartheta}(|\tilde{\psi}'|^2,\tilde{\psi}')\tilde{\psi}'
ds,\,\, \bar{d}= \int_0^1
\tilde{\Theta}(|\tilde{\psi}'|^2,\tilde{\psi}') ds,
\end{eqnarray*}
with $\tilde{\psi}=s\bar{\psi}_1+(1-s)\bar{\psi}_2$, where
$\tilde{\Theta}$ and $\tilde{\vartheta}$ are defined in
(\ref{Ederivativeeqright1}) and (\ref{Ederivativeeqright2})
respectively. Multiplying $\bar{\phi}$ on both sides of the equation
in (\ref{EblowupODEdifference}), and integrating it over $[0,1]$, we
have
\begin{equation*}
\int_0^1\frac{|\bar{\phi}'|^2}{\tilde{H}(|\tilde{\psi}'|^2,\tilde{\psi})}dx_2\leq
-\int_0^1 \left((\tilde{F}'(\tilde{\psi}))^2+
\tilde{F}(\tilde{\psi})\tilde{F}''(\tilde{\psi})\right) \tilde{H}
(|\tilde{\psi}'|^2, \tilde{\psi}) \bar{\phi}^2 dx_2.
\end{equation*}
Note that $\|\tilde{F}'\|_{C^{0,1}(\mathbb{R}^1)}\leq C
\delta^{1-3\gamma}$, thanks to the smallness of $\delta$ and the
Poincare inequality, one has
\begin{equation*}
\int_0^1 |\bar{\phi}'|^2\leq 0.
\end{equation*}
Therefore, $\bar{\phi}=0$. So the solution to
(\ref{EblowupODEBVP}) is unique. On the other hand, by the
definition of $\tilde{H}$ and $\tilde{F}$, one knows that the
boundary value problem (\ref{EblowupODEBVP}) has a solution
\begin{equation*}
\bar{\psi}=\bar{\psi}(x_2)=\int_0^{x_2}\rho_0u_0(s)ds.
\end{equation*}
This finishes the proof of the Lemma.
\end{pf}

It follows from Lemma \ref{Elemmaasymptotic} and
(\ref{Eblowup2ndHolderconvergence}) that the flow induced by the
stream function satisfies (\ref{Easymptoticupstream1}) and
(\ref{Easymptoticupstream2}).

The asymptotic behavior in the downstream can be obtained by a
similar argument.

An important direct consequence of this far fields behavior is a
better maximum estimate for the stream function.
\begin{prop}\label{Elemmastreammaxestimate}
If $\|B'\|_{C^{0,1}([0,1])}=\delta\leq \min\{\delta_1,
\delta_2\}$, $B'(0)\leq 0$, $B'(1)\geq 0$ and $m\in
(\delta^{\gamma},\bar{m})$, then the solution established in
Proposition \ref{Elemmaexistence} satisfies
(\ref{Estreambasicestimate}).
\end{prop}
\begin{pf}
It follows from Proposition \ref{Elemmaasymptotic} that
\begin{eqnarray*}
\psi(x_1,x_2)\rightarrow \int_0^{x_2}\rho_0
u_0(s)ds\,\,\text{uniformly as}
\,\, x_1\rightarrow -\infty,\\
\psi(x_1,x_2)\rightarrow \int_a^{x_2}\rho_1
u_1(s)ds\,\,\text{uniformly as} \,\, x_1\rightarrow +\infty.
\end{eqnarray*}
Therefore, for any $\epsilon>0$, there exists $L>0$ such that
\begin{equation}\label{Estreamouterestimate}
-\epsilon\leq \psi(x_1,x_2)<m+\epsilon\,\,\text{if}\,\, |x_1|\geq L.
\end{equation}
Note that $\tilde{F}'(\psi)\geq 0$ in the domain $\{\psi\geq m\}$,
thus
\begin{equation*}
\tilde{A}_{ij}(D\psi,\psi)\partial_{ij}\psi\geq  0,\text{  in the
domain } \{\psi\geq m\}\cap\{|x_1|\leq L\},
\end{equation*}
where $\tilde{A}_{ij}$ is defined in
(\ref{Esecondordercoefficients}). By maximum principle, one has
\begin{equation*}
-\epsilon\leq \psi(x_1,x_2)\leq
m+\epsilon\,\,\text{in}\,\,\{\psi\geq m\}\cap\{|x_1|\leq L\}.
\end{equation*}
Since  $\tilde{F}'(\psi)\leq 0$ in the domain $\{\psi\leq 0\}$,
thus, similarly, one can show that
\begin{equation*}
-\epsilon\leq \psi(x_1,x_2)\leq
m+\epsilon\,\,\text{in}\,\,\{\psi\leq 0\}\cap\{|x_1|\leq L\}.
\end{equation*}
Combining these estimates with (\ref{Estreamouterestimate}), it
yields
\begin{equation*}
-\epsilon\leq \psi(x_1,x_2)\leq m+\epsilon\,\,\text{in}\,\,\Omega.
\end{equation*}
Since $\epsilon$ is arbitrary, one has
\begin{equation*}
0\leq \psi(x_1,x_2)\leq m\,\,\text{in}\,\,\Omega.
\end{equation*}
This finishes the proof of the Proposition.
\end{pf}

It follows from estimates (\ref{Estreambasicestimate}) and
(\ref{Estreamgradientestimate}) that the solutions established in
Proposition \ref{Elemmaexistence} are solutions to
(\ref{Estreampb}) when the assumptions of Proposition
\ref{Elemmastreammaxestimate} are satisfied.

In fact, one can also use energy estimates to show that uniformly
subsonic solution to (\ref{Estreampb}) is unique.
\begin{prop}\label{Elemmaunique}
Let the boundary of $\Omega$ satisfy
(\ref{Enontrivialbdy})-(\ref{Eboundary3}). Then there exists
$\delta_3\in(0,\bar{\delta}_0]$ such that if
\begin{enumerate}
\item[(i).] $\|B'\|_{C^{0,1}([0,1])}=\delta\leq \delta_3$,
\item[(ii).] $m\in (\delta^{\gamma},\bar{m})$,
\end{enumerate}
then there exists at most one solution $\psi$ to (\ref{Estreampb})
satisfying
\begin{equation}\label{Euniquecondition}
0\leq\psi\leq m,\,\, |\nabla\psi|^2-\Sigma^2(\mathcal{B}(\psi)) \leq
-\epsilon \,\,\text{for some}\,\, \epsilon>0,
\end{equation}
where $H$ and $F$ are defined by $B$ and $m$ as in Section
\ref{SEreformulation}, and $
\mathcal{B}(\psi)=h(\rho_0)+\frac{F^2(\psi)}{2}$,
\end{prop}
\begin{pf}
Let $\psi_1$ and $\psi_2$ be two solutions to (\ref{Estreampb}). Set
$\psi=\psi_1-\psi_2$. Then $\psi$ satisfies
\begin{equation}\label{Euniquediffeq}
\left\{
\begin{array}{ll}
\partial_{i}(a_{ij}\partial_j\psi)+\partial_i(b_i\psi)=
c_i\partial_i\psi+d\psi,\\
\psi=0\,\, \text{on} \,\, S_1\bigcup S_2,
\end{array}
\right.
\end{equation}
where
\begin{eqnarray*}
a_{ij}&=&\int_0^1 \frac{A_{ij}(D\tilde{\psi},\tilde{\psi})}
{H^2(|\nabla\tilde{\psi}|,\tilde{\psi})} ds,\,\, b_i=\int_0^1
\frac{-H_2 (|\nabla\tilde{\psi}|^2,
\tilde{\psi})\partial_i\tilde{\psi}
}{H^2(|\nabla\tilde{\psi}|,\tilde{\psi})}
ds,\\
c_i&=&\int_0^1 \vartheta(|\nabla\tilde{\psi}|^2,
\tilde{\psi})\partial_i\tilde{\psi} ds,\,\, d= \int_0^1
\Theta(|\nabla\tilde{\psi}|^2, \tilde{\psi}) ds,
\end{eqnarray*}
here $\tilde{\psi}=s\psi_1+(1-s)\psi_2$, $A_{ij}$, $\Theta$ and
$\vartheta$ are defined similar to (\ref{Eprincipalcoefficient}),
(\ref{Ederivativeeqright1}) and (\ref{Ederivativeeqright2}),
respectively except we replace $\tilde{F}$ and  $\tilde{H}$ by $F$
and $H$.

Multiplying $\eta^2\psi^{+}$ on both sides of equation in
(\ref{Euniquediffeq}), where $\eta$ is defined in
(\ref{Ecutofftestfunction})  and $\psi^{+}(x)=\max\{\psi(x),0\}$,
then similar to the proof of Lemma \ref{Elemmaasymptotic}, one has
\begin{equation*}
\iint_{\Omega\cap\{|x_1|\leq l\}\cap \{\psi\geq
0\}}|\nabla\psi|^2dx_1dx_2\leq C(\underline{B},\epsilon)
\iint_{\Omega\cap\{l\leq |x_1|\leq l+1\}\cap \{\psi\geq
0\}}|\nabla\psi|^2dx_1dx_2.
\end{equation*}
It follows from Lemma \ref{Elemmaasymptotic} that $\psi_1$ and
$\psi_2$ have the same far fields behavior, thus $|\psi|$ and
$|\nabla\psi|\rightarrow 0$ as $|x_1|\rightarrow \infty$. Thus
\begin{equation*}
\iint_{\Omega\cap\{\psi\geq 0\}}|\nabla\psi|^2=0,
\end{equation*}
so $\psi\leq 0$. Similarly, one can show that $\psi\geq 0$.
Therefore, $\psi=0$. This finishes the proof of the Proposition.
\end{pf}

\section{Refined Estimates for the Boundary Value Problem
for Stream Functions}\label{SErefined}

In this section, we will derive some refined estimates for solutions
to the problem (\ref{Estreampb}). Combining these refined estimates
with the estimates obtained in section \ref{SEexistence} and section
\ref{SEasymptotic}, one will get a solution to the Euler equations
(\ref{EEulercontinuityeq})-(\ref{EEulermomentumeq2}), with the
boundary condition (\ref{Enoflowbc}) and the constrains
(\ref{Emassflux}) and (\ref{EasymtoticBernoulli}). More precisely,
it will be shown that $\psi_{x_2}$ is always positive, therefore,
$u=\psi_{x_2}/\rho=\psi_{x_2}/H(|\nabla\psi|^2,\psi)$ satisfies
(\ref{Esubsonic}). This positivity of the horizontal velocity and
the asymptotic behavior yield that $(\rho, u,v)$ induced by $\psi$
satisfies the original Euler equations, the  boundary conditions and
the constrains on mass flux and Bernoulli's constant.

\begin{lem}\label{Elemmarefined}
Let the boundary of $\Omega$ satisfy
(\ref{Enontrivialbdy})-(\ref{Eboundary3}). Then there exists
$\delta_4\in (0,\bar{\delta}_0]$ such that if
\begin{enumerate}
\item[(i).] $\|B'\|_{C^{0,1}([0,1])}=\delta\leq \delta_4$,
\item[(ii).] $m\in (\delta^{\gamma},\bar{m})$, \item[(iii).] $\psi$
satisfies (\ref{Euniquecondition}) and solves (\ref{Estreampb})
with $H$ and $F$ defined  by $B$ and $m$ as in Section
\ref{SEreformulation},
\end{enumerate}
then $\psi$ satisfies
\begin{equation}\label{Estreambasicestimateagain}
0< \psi< m\,\,\text{in}\,\,\Omega,
\end{equation}
and
\begin{equation}
\psi_{x_2}>0\,\,\text{in}\,\,\bar{\Omega}.
\end{equation}
\end{lem}

\begin{pf}
It follows from (\ref{Euniquecondition}) that $\psi$ achieves its
minimum on $S_1$ and maximum on $S_2$, hence
\begin{equation}\label{Ehvelocityroughbdy}
\psi_{x_2}\geq 0\,\,\text{on}\,\,\partial\Omega.
\end{equation}
On the other hand, $U=\psi_{x_2}$ satisfies
\begin{equation}\label{Ehorizontalvelocityeq}
\begin{array}{ll}
\,\,\,\,\,\,\partial_{i}\left(\frac{A_{ij}(D\psi,\psi)}{H^2(|\nabla\psi|^2,\psi)}\partial_j
U\right)-\partial_i\left(\frac{H_2(|\nabla\psi|^2,\psi)\partial_i\psi}{H^2(|\nabla\psi|^2,\psi)}U\right)\\
=\Theta(|\nabla\psi|^2,\psi)U
+\vartheta(|\nabla\psi|^2,\psi)\partial_i\psi\partial_i U
\end{array}
\end{equation}
in the weak sense, where $A_{ij}$, $\Theta$ and $\vartheta$ are
defined similar to (\ref{Eprincipalcoefficient}),
(\ref{Ederivativeeqright1}) and (\ref{Ederivativeeqright2})
respectively except we replace $\tilde{F}$ and  $\tilde{H}$ by $F$
and $H$. We first claim that
\begin{equation}\label{Ehvelocityroughinterior}
U\geq 0 \,\,\text{in}\,\,\Omega.
\end{equation}

Indeed, it follows from Lemma \ref{Elemmaasymptotic} that
$U(x_1,x_2)>0$ when $|x_1|>l$ for some $l$ sufficiently large.
Multiplying (\ref{Ehorizontalvelocityeq}) by $U^{-}=\min\{U,0\}$,
and using (\ref{Ehvelocityroughbdy}), one may get that
\begin{eqnarray*}
&&\iint_{\{U\leq 0\}} \frac{|\nabla U|^2}{H(|\nabla\psi|^2,\psi)}dx_1 dx_2\\
&=&\iint_{\{U\leq 0\}} \frac{2H_1(|\nabla\psi|^2,\psi)}{H^2
(|\nabla\psi|^2,\psi)}|\nabla\psi\cdot\nabla
U|^2 dx_1dx_2\\
&&+ \iint_{\{U\leq 0\}}
\frac{H_2(|\nabla\psi|^2,\psi)\partial_i\psi}
{H^2(|\nabla\psi|^2,\psi)} U\partial_i U dx_1dx_2\\
&&-\iint_{\{U\leq 0\}}  \Theta(|\nabla\psi|^2,\psi)
U^2dx_1dx_2-\iint_D \vartheta(|\nabla\psi|^2,\psi)
U\nabla\psi\cdot\nabla U dx_1dx_2\\
&\leq &-\iint_{\{U\leq 0\}} (F^{''}(\psi)F(\psi)+(F'(\psi))^2)
H(|\nabla\psi|^2,\psi)U^2dx_1dx_2\\
&\leq & C\delta^{1-3\gamma}\iint_{\{U< 0\}} U^2 dx_1dx_2.
\end{eqnarray*}
Define $K_{x_1}=\{x_2|f_1(x_1)\leq x_2\leq f_2(x_1),
U(x_1,x_2)<0\}$, then $K_{x_1}$ is an open set for each $x_1$. Let
$K_{x_1}=\cup_{i\in \mathscr{S}}I_{x_1}^i$, where $I_{x_1}^i$ are
connected components of $K_{x_1}$. For each $x_2\in I_{x_1}^i$,
\begin{equation*}
U(x_1,x_2)=\int_{\min I_{x_1}^i}^{x_2}U(x_1,s)ds.
\end{equation*}
Therefore,
\begin{equation*}
\begin{aligned}
&\iint_{\{U<0\}}U^2(x_1,x_2)dx_1dx_2\\
=&\int_{-l}^ldx_1\Sigma_{i\in\mathscr{S}}
\int_{I_{x_1}^i}U^2(x_1,x_2)dx_2\\
=&\int_{-l}^ldx_1\Sigma_{i\in\mathscr{S}}
\int_{I_{x_1}^i}\left(\int_{\min
I_{x_1}^i}^{x_2}\partial_{x_2}U(x_1,s)ds\right)^2dx_2\\
\leq& \int_{-l}^ldx_1\Sigma_{i\in\mathscr{S}}
\int_{I_{x_1}^i}\int_{\min I_{x_1}^i}^{\max
I_{x_1}^i}(\partial_{x_2}U(x_1,s))^2ds(\max I_{x_1}^i-\min
I_{x_1}^i)dx_2\\
=& \int_{-l}^ldx_1\Sigma_{i\in\mathscr{S}}(\max I_{x_1}^i-\min
I_{x_1}^i)^2 \int_{\min I_{x_1}^i}^{\max
I_{x_1}^i}(\partial_{x_2}U(x_1,s))^2ds\\
\leq&
\max_{x_1\in\mathbb{R}}|f_2(x_1)-f_1(x_1)|^2\int_{-l}^ldx_1\Sigma_{i\in\mathscr{S}}
\int_{\min I_{x_1}^i}^{\max I_{x_1}^i}(\partial_{x_2}U(x_1,s))^2ds\\
\leq& \max_{x_1\in\mathbb{R}}|f_2(x_1)-f_1(x_1)|^2 \iint_{\{U< 0\}}
|\nabla U|^2 dx_1dx_2
\end{aligned}
\end{equation*}
Hence,
\begin{equation*}
\iint_{\{U\leq 0\}} \frac{|\nabla U|^2}{H(|\nabla\psi|^2,\psi)}dx_1
dx_2 \leq C\delta^{1-3\gamma}\iint_{\{U\leq 0\}} |\nabla U|^2
dx_1dx_2,
\end{equation*}
which implies
\begin{equation*}
\iint_{\{U\leq 0\}}|\nabla U|^2 dx_1dx_2\leq 0,
\end{equation*}
so (\ref{Ehvelocityroughinterior}) must hold.

Now, we use an argument similar to the proof of Lemma 1 in \S9.5.2
in \cite{Evans} to show that
\begin{equation}\label{Eestimatenostagnation}
\psi_{x_2}=U>0\,\,\text{in}\,\,\Omega
\end{equation}
holds for any weak solutions $U$ to (\ref{Ehorizontalvelocityeq}).

Indeed, let $\tilde{U}=e^{-\sigma x_2}U$. Then $\tilde{U}$ is a
nonnegative weak solution to
\begin{equation*}
\partial_i\left(\frac{A_{ij}}{H^2}e^{\sigma
x_2}\partial_j\tilde{U}\right)+\left( \frac{A_{i2}}{H^2}\sigma-
\frac{H_2\partial_i\psi}{H^2}-\vartheta(|\nabla\psi|^2,\psi)
\partial_i\psi\right)e^{\sigma x_2}\partial_i \tilde{U}+G e^{\sigma x_2}\tilde{U}=0,
\end{equation*}
where $A_{ij}$ and $\vartheta$ are defined in
(\ref{Eprincipalcoefficient}) and (\ref{Ederivativeeqright2}), and
\begin{equation*}
G=
\frac{A_{22}}{H^2}\sigma^2+\left(\partial_i\left(\frac{A_{i2}}{H^2}\right)
-\frac{H_2\partial_2\psi}{H^2}
-\vartheta(D\psi,\psi)\partial_2\psi\right)\sigma\\
-\partial_i\left(\frac{H_2\partial_i\psi}{H^2}\right)
-\Theta(|\nabla\psi|^2,\psi)
\end{equation*}
with $\Theta$ defined in (\ref{Ederivativeeqright1}). Choose
$\sigma>0$ sufficiently large. Then $G>0$. Thus
\begin{equation*}
\partial_i\left(\frac{A_{ij}}{H^2}e^{\sigma
x_2}\partial_j\tilde{U}\right)+\left( \frac{A_{i2}}{H^2}\sigma-
\frac{H_2\partial_i\psi}{H^2}-\vartheta(|\nabla\psi|^2,\psi)
\partial_i\psi\right)e^{\sigma x_2}\partial_i \tilde{U}\leq 0.
\end{equation*}
It follows from Theorem 8.19 in \cite{GT} that
(\ref{Eestimatenostagnation}) holds.

Now, (\ref{Estreambasicestimateagain}) follows directly from
(\ref{Euniquecondition}) and (\ref{Eestimatenostagnation}).

Since $\psi=m$ on $S_2$, if $F'(m)> 0$, then for any
$(x_1^0,f_2(x_1^0))\in S_2$, there exists a small disk
$\mathcal{N}\subset\Omega$ satisfying $\bar{\mathcal{N}}\cap
\bar{\Omega}=(x_1^0,f_2(x_1^0))$ such that $F'(\psi)\geq 0$ in
$\mathcal{N}$, therefore,
\begin{equation*}
A_{ij}(D\psi,\psi)\partial_{ij}\psi\geq
0\,\,\text{in}\,\,\mathcal{N}.
\end{equation*}
Moreover, by (\ref{Estreambasicestimateagain}), $\psi< m$ in
$\mathcal{N}$. Thus, by the Hopf Lemma, one has
$\psi_{x_2}(x_1^0,f_2(x_1^0))>0$.

In the case $F'(m)=0$, $\psi$ satisfies
\begin{equation*}
A_{ij}(D\psi,\psi)\partial_{ij}(\psi-m)+R(\psi-m)=0,
\end{equation*}
where
\begin{equation*}
R=-\frac{F(\psi)H^5c^2} {H^2(|\nabla\psi|^2,\psi)c^2
-|\nabla\psi|^2} \frac{F'(\psi)-F'(m)}{\psi-m}.
\end{equation*}
It follows from the Hopf lemma that
\begin{equation*}
\partial_{x_2}\psi>0\,\,\text{in}\,\, S_2.
\end{equation*}

Similarly, one can show that $\psi_{x_2}(x_1, f_1(x_1))>0$ for any
$x_1\in\mathbb{R}$.

This finishes the proof of the Lemma.
\end{pf}

Choose $\delta_0=\min\{\delta_1,\delta_2,\delta_3,\delta_4\}$,
then $\delta_0>0$. If $\|B'\|_{C^{0,1}([0,1])}=\delta\leq
\delta_0$, for any $m\in (\delta^{\gamma},2\delta_0^{\gamma/2})$,
there exists a solution to the problem (\ref{Estreampb}). It
follows from Lemma \ref{Elemmaasymptotic} and Lemma
\ref{Elemmarefined} that the flow field induced by $\psi$
satisfies (\ref{Ephorizontalvelocity}) and
(\ref{EcequivstreamEuler}), hence Proposition
\ref{Elemmaequivalent} guarantees the existence of Euler flows.
Furthermore, Proposition \ref{Elemmaequivalent} and Proposition
\ref{Elemmaunique} imply uniqueness of Euler flows with asymptotic
condition (\ref{EasymtoticBernoulli}), mass flux condition
(\ref{Emassflux}), (\ref{Esubsonic}), and asymptotic behavior
(\ref{Easymptoticupstream1})-(\ref{Easymptoticupstream2}).

\section{Existence of Critical Mass Flux}\label{SEcritical}
So far, we have shown that, for the given Bernoulli's function in
the upstream satisfying (\ref{EassumptiononBernoulli}), there exist
Euler flows as long as $m\in
(\delta^{\gamma},2\delta_0^{\gamma/2})$. In this section, we will
increase $m$ as large as possible.

\begin{prop}\label{Epropcritical}
Let $\Omega$ satisfy (\ref{Enontrivialbdy})-(\ref{Eboundary3}) and
$B$ satisfy (\ref{EassumptiononBernoulli}) and
(\ref{Ecricondition}). Then there exists $\hat{m}\leq \bar{m}$ such
that if $m\in (\delta^{\gamma},\hat{m})$, there exists a unique
$\psi$ which satisfies
\begin{equation}\label{Eellipticity}
0< \psi< m\,\,\text{in}\,\,\Omega,\,\,\text{and}\,\,
M(m)=\sup_{\bar{\Omega}}\left(|\nabla\psi|^2-\Sigma^2(\mathcal{B}(\psi))\right)<0,
\end{equation}
and solves (\ref{Estreampb}), where
$\mathcal{B}(\psi)=h(\rho_0)+F^2(\psi)/2$. Furthermore, either
$M(m)\rightarrow 0$ as $m\rightarrow \hat{m}$, or there does not
exist $\sigma>0$ such that (\ref{Estreampb}) has solutions for all
$m\in (\hat{m},\hat{m}+\sigma)$ and
\begin{equation}\label{Ebifurcationestimate}
\sup_{m\in(\hat{m},\hat{m}+\sigma)}M(m)<0.
\end{equation}
\end{prop}
\begin{pf}
The basic idea of the proof for Proposition is quite similar to
that in \cite{Bers1,XX1}.

For the given Bernoulli's function $B$ in the upstream satisfying
(\ref{EassumptiononBernoulli}) and any $m\in
(\delta^{\gamma},\bar{m})$, one can define $\rho_0$ and $u_0(x_2)$,
and therefore $F(\psi)$ as in Section \ref{SEreformulation}. Note
that $\rho_0$ and $F$ depend on $m$ by definition, thus in this
section, we will denote them by $\rho_0(m)$ and $F(\psi;m)$
respectively.

When $B$ satisfies (\ref{Ecricondition}), one has
\begin{equation*}
F'(m)=F'(0)=0.
\end{equation*}
Thus $\tilde{F}'$, the extension of $F'$ is Section 3, has the
following simple form
\begin{equation}
\tilde{F}'(s)=\left\{
\begin{array}{ll}
F'(s),\qquad \text{if}\,\, 0\leq s\leq m,\\
0,\qquad \text{if}\,\, s<0\,\,\text{or}\,\, s>m.
\end{array}
\right.
\end{equation}
Set $\tilde{F}(s)=\int_0^s \tilde{F}'(s)ds$. Then it is easy to
check that
\begin{equation}\label{EspecialestF}
B_0<\underline{B}\leq h(\rho_0)+\frac{\tilde{F}^2(s)}{2}\leq
\bar{B}\,\, \text{and}\,\,
\|\tilde{F}'\|_{C^{0,1}(\mathbb{R}^1)}\leq C\delta^{1-2\gamma}.
\end{equation}

Let $\{\varepsilon_n\}_{n=1}^{\infty}$ be a strictly decreasing
sequence of positive numbers such that $\varepsilon_1\leq
\varepsilon_0/4$ and $\varepsilon_n\downarrow 0$. One can truncate
$H$ associated with $\varepsilon_n$ as follows
\begin{equation}
\tilde{H}^{(n)}(|\nabla\psi|^2,\psi;m)=J(\tilde{\Delta}_n(|\nabla\psi|^2,\psi;m),
\tilde{\mathcal{B}}_n(\psi;m)).
\end{equation}
To give a clear explanation of this definition, we first introduce
a sequence of smooth increasing functions $\zeta_n$ such that
\begin{equation*}
\zeta_n(s)=\left\{
\begin{array}{ll}
s,\,\, &\text{if}\,\, s<-2\varepsilon_n,\\
-\varepsilon_n,\,\,&\text{if}\,\,s\geq -\varepsilon_n.
\end{array}
\right.
\end{equation*}
Now one can define
\begin{equation*}
\tilde{\Delta}_n
(|\nabla\psi|^2,\psi;m)=\zeta_n(|\nabla\psi|^2-\Sigma^2(\tilde{\mathcal{B}}_n(\psi;m)))
+\Sigma^2(\tilde{\mathcal{B}}_n(\psi;m)),
\end{equation*}
where
\begin{equation*}
\tilde{\mathcal{B}}_n(\psi;m)=h(\rho_0(m))+\frac{\tilde{F}^2(\psi;m)}{2}.
\end{equation*}

It is easy to see that there exist two positive constants
$\lambda(n)$ and $\Lambda(n)$ such that
\begin{equation*}
\lambda(n)|\xi|^2\leq \tilde{A}^{(n)}_{ij}(q,z;m) \xi_i\xi_j \leq
\Lambda(n)|\xi|^2
\end{equation*}
for any $z\in\mathbb{R}^1$, $q\in \mathbb{R}^2$ and
$\xi\in\mathbb{R}^2$, where
\begin{equation*}
\tilde{A}^{(n)}_{ij}(q,z;m)=\tilde{H}^{(n)}(|q|^2,z;m)\delta_{ij}-
2\tilde{H}^{(n)}_1(|q|^2,z;m)q_i q_j.
\end{equation*}
Thus it follows from the argument in Section \ref{SEexistence}
that  for any $m\in(\delta^{\gamma},\bar{m})$, there exists a
solution $\psi^{(n)}(x;m)$ to the problem
\begin{equation}\label{Ecriticalcutproblem}
\left\{
\begin{array}{ll}
\tilde{A}_{ij}^{(n)}(D\psi,\psi;m)\partial_{ij}\psi=\tilde{\mathcal{F}}_n(D\psi,\psi;m)
\,\,\mathrm{in}\,\,\Omega,\\
\psi=\frac{x_2-f_1(x_1)}{f_2(x_1)-f_1(x_1)}m\,\,\text{on}\,\,\partial\Omega,
\end{array}
\right.
\end{equation}
where
\begin{equation*}
\tilde{\mathcal{F}}_n=\tilde{F}\tilde{F}' \tilde{H}^{(n)}
\left(\frac{ ((\tilde{H}^{(n)})^2+\Sigma\Sigma'(\zeta_n'-1))
\tilde{\Delta}_n} {(\tilde{H}^{(n)})^2 c^2 -\tilde{\Delta}_n}
+(\tilde{H}^{(n)})^2\right),
\end{equation*}
where we ignore some obvious independent variables in the
definition of $\tilde{\mathcal{F}}_n$. Moreover, if
\begin{equation}\label{Ecriticaluniformellipticity}
|\nabla \psi^{(n)}|^2-\tilde{\mathcal{B}}(\psi^{(n)};m)\leq
-2\varepsilon_n,
\end{equation}
then $\zeta_n'=1$. Similar to Section \ref{SEexistence}, one has
\begin{equation*}
0\leq \psi^{(n)}(x;m)\leq m.
\end{equation*}

Since $\tilde{F}$ satisfies (\ref{EspecialestF}) independent of
$\varepsilon_n$, one can estimate  $I_5$  in (\ref{Ecptformkeyest})
as that in (\ref{EestI5}). Furthermore, it follows from the same
arguments in Lemma \ref{Elemmaasymptotic} that the solution to
(\ref{Ecriticalcutproblem}) satisfying
(\ref{Ecriticaluniformellipticity}) has far fields behavior as
(\ref{Eblowupexplicitform}). In addition, by Proposition
\ref{Elemmaunique}, such a solution is unique among the class of
solutions satisfying (\ref{Eblowupexplicitform}).

Note that in general, we do not know uniqueness of solutions to
problem (\ref{Ecriticalcutproblem}). Set
\begin{equation}
S_n(m)=\{\psi^{(n)}(x;m)|\psi^{(n)}(x;m)\,\,\text{solves the
problem}\,\,(\ref{Ecriticalcutproblem})  \}.
\end{equation}
Define
\begin{equation}
M_n(m)=\inf_{\psi^{(n)}\in S_n(m)}\sup_{\bar{\Omega}}
(|\nabla\psi^{(n)}(x;m)|^2-\Sigma^2(\tilde{\mathcal{B}}(\psi^{(n)};m))),
\end{equation}
and
\begin{equation*}
T_n=\{s|\delta_0^{\gamma}\leq s\leq \bar{m}, M_n(m)\leq
-4\varepsilon_n\,\,\text{if}\,\,m\in (\delta^{\gamma},s)\}.
\end{equation*}
It follows from Proposition \ref{Elemmaexistence}, Lemma
\ref{Elemmaasymptotic} and Proposition \ref{Elemmastreammaxestimate}
that $[\delta_0^{\gamma},2\delta_0^{\gamma/2}]\subset T_n$,
therefore, $T_n$ is not an empty set. Define $m_n=\sup T_n$.

The sequence $\{m_n\}$ has some nice properties.

First, $M_n(m)$ is left continuous for $m\in (\delta^{\gamma},m_n]$.
Indeed, let $\{m_n^{(k)}\}\subset(\delta^{\gamma}, m_n)$ and
$m_n^{(k)}\uparrow m$. Since $M_n(m_n^{(k)})\leq -4\varepsilon_n$,
one has
\begin{equation*}
\|\psi^{(n)}(x;m_n^{(k)})\|_{C^{2,\alpha}(\bar{\Omega})}\leq C.
\end{equation*}
Therefore, there exists a subsequence $\psi^{(n)}(x;m_n^{(k_l)})$
such that $\psi^{(n)}(x;m_n^{(k_l)})\rightarrow \psi$, moreover,
$\psi$ solves (\ref{Ecriticalcutproblem}). Thus $M_n(m)\leq \lim
M_n(m_n^{(k_l)})$. So $M_n(m)\leq -4\varepsilon_n$. Note that all
these solutions satisfy the far fields behavior as
(\ref{Eblowupexplicitform}), by uniqueness of solutions in this
class, $M_n(m)= \lim M_n(m_n^{(k)})$.

Second, $m_n< \bar{m}$. Suppose not, by the definition of $m_n$,
$\bar{m}\in T_n$. It follows from the left continuity of $M_n$,
$M_n(\bar{m})\leq -4\varepsilon_n$. Thus by means of the proof of
Lemma \ref{Elemmaasymptotic}, $\psi^{(n)}(x;\bar{m})$ has far field
behavior as in (\ref{Eblowupexplicitform}). However, it follows from
the definition of $\bar{m}$ that
\begin{eqnarray*}
& &\sup_{x\in\bar{\Omega}}\left( |\nabla\psi^{(n)}(x;\bar{m})|^2-
\Sigma^2(\tilde{\mathcal{B}}_n(\psi^{(n)}(x;\bar{m})))\right)\\
&\geq&\sup_{x_2\in[0,1])}\max\{(|\rho_0(\bar{m})u_0(x_2;\bar{m})|^2-\Sigma^2(B(x_2))),
(|\rho_1(\bar{m})u_1(y(x_2);\bar{m})|^2-\Sigma^2(B(x_2)))\}\\
&=&0,
\end{eqnarray*}
where $y=y(s)$ is the function defined in (\ref{Eflowmap}). Thus
$M_n(\bar{m})\geq 0$. This is a contradiction. Therefore $m_n<
\bar{m}$.

Finally, $\{m_n\}$ is an increasing sequence. This follows from
the definition of $\{m_n\}$ directly.

Define $\hat{m}=\lim_{n\rightarrow \infty} m_n$. Based on previous
properties of $\{m_n\}$, $\hat{m}$ is well-defined and
$\hat{m}\leq \bar{m}$.

Note that for any $m\in (\delta^{\gamma},\hat{m})$, there exists
$m_n>m$, therefore $M_n(m)\leq -4\varepsilon_n$. Thus
$\psi=\psi^{(n)}(x;m)$ solves (\ref{Estreampb}) and
\begin{equation*}
\sup_{\bar{\Omega}}
(|\nabla\psi|^2-\Sigma^2(\mathcal{B}(\psi)))=M_n(m)\leq
-4\varepsilon_n.
\end{equation*}

If $\sup_{m\in(\delta^{\gamma},\hat{m})}M(m)<0$, then there exists
$n$ such that
$\sup_{m\in(\delta^{\gamma},\hat{m})}M(m)<-4\varepsilon_n$. As the
same as the proof for the left continuity of $M_n(m)$ on
$(\delta^{\gamma},m_n]$, $M_n(\hat{m})\leq-4\varepsilon$. Suppose
that there exists $\sigma>0$ such that (\ref{Estreampb}) always has
a solution $\psi$ for $m\in(\hat{m},\hat{m}+\sigma)$, and
\begin{equation}
\sup_{m\in(\hat{m},\hat{m}+\sigma)}
M(m)=\sup_{m\in(\hat{m},\hat{m}+\sigma)}\sup_{\bar{\Omega}}
(|\nabla\psi|^2-\Sigma^2(\mathcal{B}(\psi)))< 0.
\end{equation}
Then there exists $k>0$ such that
\begin{equation*}
\sup_{m\in(\hat{m},\hat{m}+\sigma)}
M(m)=\sup_{m\in(\hat{m},\hat{m}+\sigma)}\sup_{\bar{\Omega}}
(|\nabla\psi|^2-\Sigma^2(\mathcal{B}(\psi)))\leq
-4\varepsilon_{n+k}.
\end{equation*}
This yields that $m_{n+k}\geq \hat{m}+\sigma$. So there is a
contradiction. The contradiction implies that either
$M(m)\rightarrow 0$, or there does not exist $\sigma>0$ such that
(\ref{Estreampb}) has solution for all $m\in
(\hat{m},\hat{m}+\sigma)$ and (\ref{Ebifurcationestimate}) holds.

This finishes the proof of the Proposition.
\end{pf}

It follows from Lemma \ref{Elemmaasymptotic}, Lemma
\ref{Elemmarefined} and Proposition \ref{Epropcritical} that if $B$
satisfies (\ref{EassumptiononBernoulli}) and
$m\in(\delta^{\gamma},\hat{m})$, then there exists an Euler flow
through the nozzle. Collecting all results obtained together, we
complete the proof of Theorem \ref{EThexistence}.

\bigskip
{\bf Acknowledgement.} Part of this work was done when the first
author was visiting the Center for Nonlinear Studies in Northwestern
University, Xi'an, China, he would like to thank the center's
hospitality and support during his visit. The authors would like to
thank Professor Tao Luo and Mr. Wei Yan for helpful discussions.
This research is supported in part by Hong Kong RGC Earmarked
Research Grants CUHK04-28/04, CUHK-4040/06, CUHK-4042/08P and RGC
Central Allocation Grant CA05/06.SC01.


\bigskip

\begin{thebibliography}{99}



\bibitem{Alber}
H. D. Alber, {\it Existence of three-dimensional, steady, inviscid,
incompressible flows with nonvanishing vorticity}, Math. Ann.,
292(1992), no. 3, 493--528.



\bibitem{Bers1}
L. Bers, {\it Existence and uniqueness of a subsonic flow past a
given profile}, Comm. Pure Appl. Math.,  7(1954), 441--504.



\bibitem{Bers2}
L. Bers, {\it Mathematical aspects of subsonic and transonic gas
dynamics}, John Wiley \& Sons, Inc., 1958.


\bibitem{CCS}
Gui-Qiang Chen, Jun Chen and  Kyungwoo Song, {\it Transonic nozzle
flows and free boundary problems for the full Euler equations}, J.
Differential Equations, 229(2006), no. 1, 92--120.



\bibitem{CDSW}
Gui-Qiang Chen, Constantine M. Dafermos, Marshall Slemrod and Dehua
Wang,  {\it On two-dimensional sonic-subsonic flow}, Comm. Math.
Phys.,  271(2007), no. 3, 635--647.



\bibitem{CF1}
Gui-Qiang Chen and Mikhail Feldman,  {\it Multidimensional transonic
shocks and free boundary problems for nonlinear equations of mixed
type}. J. Amer. Math. Soc., 16(2003), no. 3, 461--494.

\bibitem{CF2}
Gui-Qiang Chen and Mikhail Feldman,  {\it Steady transonic shocks
and free boundary problems for the Euler equations in infinite
cylinders}. Comm. Pure Appl. Math.,  57(2004), no. 3, 310--356.


\bibitem{ChenJun}
Jun Chen,  {\it Subsonic Euler Flows in Half Plane}, preprint, 2007.


\bibitem{CF}
R. Courant and K. O. Friedrichs, {\it Supersonic flow and shock
waves},  Applied Mathematical Sciences, Vol. 21. Springer-Verlag,
New York-Heidelberg, 1976.

\bibitem{Dong}
Guang Chang Dong, {\it Nonlinear partial differential equations of
second order},  Translations of Mathematical Monographs, 95,
American Mathematical Society, Providence, RI, 1991.



\bibitem{Evans}
Lawrence C. Evans, {\it Partial differential equations}, Graduate
Studies in Mathematics, 19, American Mathematical Society,
Providence, RI, 1998.




\bibitem{Feistauer}
M. Feistauer, {\it Mathematical methods in fluid dynamics}, Pitman
Monographs and Surveys in Pure and Applied Mathematics, 67, Longman
Scientific \& Technical, Harlow; copublished in the United States
with John Wiley \& Sons, Inc., New York, 1993.



\bibitem{FG1}
Robert Finn and David Gilbarg, {\it Asymptotic behavior and
uniquenes of plane subsonic flows}, Comm. Pure Appl. Math.,
10(1957), 23--63.



\bibitem{FG2}
Robert Finn and David Gilbarg, {\it Three-dimensional subsonic
flows, and asymptotic estimates for elliptic partial differential
equations}, Acta Math., 98(1957), 265--296.

\bibitem{G}
David Gilbarg, {\it Comparison methods in the theory of subsonic
flows}, J. Rational Mech. Anal.,  2(1953), 233--251.



\bibitem{GS}
David Gilbarg and M. Shiffman, {\it On bodies achieving extreme
values of the critical Mach number,  I.},  J. Rational Mech. Anal.
3(1954), 209--230.



\bibitem{GT}
David Gilbarg and Neil S. Trudinger, {\it Elliptic partial
differential equations of second order}, Second edition,
Springer-Verlag, Berlin, 1983.



\bibitem{Glass}
Olivier Glass, {\it Existence of solutions for the two-dimensional
stationary {E}uler system for ideal fluids with arbitrary force},
Ann. Inst. H. Poincar\'e Anal. Non Lin\'eaire, 20(2003), no. 6,
921--946.


\bibitem{LXY1}
Jun Li, Zhouping Xin and Huicheng Yin, {\it On transonic shocks in a
nozzle with variable end pressures}, to appear in  Communications in
Mathematical Physics, 2007.


\bibitem{LXY2}
Jun Li, Zhouping Xin and Huicheng Yin, {\it A free boundary value
problem for the Euler system and 2-D transonic shock in a large
variable nozzle}, to appear in  Mathematical Research Letters, 2009.


\bibitem{Morawetz1}
Cathleen S. Morawetz,  {\it On the non-existence of continuous
transonic flows past profiles. I.}, Comm. Pure Appl. Math., 9(1956),
45--68.



\bibitem{Morawetz2}
Cathleen S. Morawetz, {\it On the non-existence of continuous
transonic flows past profiles. II.},  Comm. Pure Appl. Math.,  10
(1957), 107--131.




\bibitem{Morawetz3}
Cathleen S. Morawetz, {\it On the non-existence of continuous
transonic flows past profiles. III.}, Comm. Pure Appl. Math.,
11(1958), 129--144.


\bibitem{Morawetz4}
Cathleen S. Morawetz, {\it Non-existence of transonic flow past a
profile}, Comm. Pure Appl. Math., 17(1964), 357--367.



\bibitem{Morawetzcpt1}
Cathleen S. Morawetz, {\it On a weak solution for a transonic flow
problem}, Comm. Pure Appl. Math., 38(1985), no. 6, 797--817.

\bibitem{Morawetzcpt2}
Cathleen S. Morawetz, {\it On steady transonic flow by compensated
compactness}, Methods Appl. Anal., 2(1995), no. 3, 257--268.



\bibitem{Rauch}
Jeffrey Rauch,  {\it B{V} estimates fail for most quasilinear
hyperbolic systems in dimensions greater than one}, Comm. Math.
Phys., 106 (1986), no. 3, 481--484.


\bibitem{Shiffman}
Max Shiffman, {\it On the existence of subsonic flows of a
compressible fluid}, J. Rational Mech. Anal., 1(1952), 605--652.

\bibitem{Sideris}
Thomas C. Sideris, {\it Formation of singularities in
three-dimensional compressible fluids}, Comm. Math. Phys.,
101(1985), no. 4, 475--485.

\bibitem{Troshkin}
O. V. Troshkin, {\it Nontraditional methods in mathematical
hydrodynamics}, Translations of Mathematical Monographs, 144,
American Mathematical Society, Providence, RI, 1995.




\bibitem{XX1}
Chunjing Xie and  Zhouping Xin, {\it Global subsonic and
subsonic-sonic flows through infinitely long nozzles}, Indiana Univ.
Math. J., 56(2007), no. 6, 2991--3023.



\bibitem{XX2}
Chunjing Xie and  Zhouping Xin, {\it Global subsonic and
subsonic-sonic flows through infinitely long axially symmetric
nozzles}, submitted to { Journal of Differential Equations}, 2007.



\bibitem{XYY}
Zhouping Xin, Wei Yan and  Huicheng Yin, {\it Transonic shock
problem for Euler system in a nozzle}, to appear in {Archive for
Rational Mechanics and Analysis}, 2008.


\bibitem{XY}
Zhouping Xin and Huicheng Yin, {\it Transonic shock in a nozzle. I.
Two-dimensional case}, Comm. Pure Appl. Math.,  58(2005), no. 8,
999--1050.


\bibitem{XY2}
Zhouping Xin and Huicheng Yin, {\it Three-dimensional transonic
shocks in a nozzle}, Pacific J. Math., 236(2008), no. 1, 139--193.



\bibitem{XY3}
Zhouping Xin and Huicheng Yin, {\it The transonic shock in a nozzle,
2-D and 3-D complete Euler systems}, J. Differential Equations,
245(2008), no. 4, 1014--1085.





\end{thebibliography}
\end{document}